\documentclass[aos,MSNbibl,seceqn,nameyear,noautosecdot,dvips]{arximspdf}

\doi{10.1214/12-AOS1077} %kopijuoti is PTS
\volume{41}
\issue{2}
\pubyear{2013}
\firstpage{802}
\lastpage{837}

\makeatletter

\newproclaim{definition}{Definition}

\newcommand{\rrvert}{\vert}
\newcommand{\llvert}{\vert}

\newtheorem{theorem}{Theorem}[section]
\newtheorem{lem}{Lemma}[section]
\newtheorem{cor}{Corollary}[section]
\newtheorem{prop}{Proposition}[section]

\newcommand{\reals}{\mathbb{R}}

\newcommand{\bbR}{\mathbb{R}}
\newcommand{\cL}{\mathcal{L}}
\newcommand{\cM}{{\mathcal{M}}}
\newcommand{\Mall}{{\mathcal{M}_{\mathrm{all}}}}
\newcommand{\tcM}{{\tilde{\mathcal{M}}}}
\newcommand{\cN}{{\mathcal{N}}}

\newcommand{\0}{{\mathbf{0}}}
\newcommand{\1}{{\mathbf{1}}}
\renewcommand{\l}{{\mathbf{l}}}
\newcommand{\ls}{{\bar{\mathbf{l}}}}
\newcommand{\A}{{\mathbf{A}}}
\newcommand{\B}{{\mathbf{B}}}
\newcommand{\D}{{\mathbf{D}}}
\newcommand{\E}{{\mathbf{E}}}
\newcommand{\F}{{\mathrm{F}}}
\newcommand{\I}{{\mathbf{I}}}
\renewcommand{\P}{{\mathbf{P}}}
\newcommand{\Q}{{\mathbf{Q}}}
\newcommand{\R}{{\mathbf{R}}}
\newcommand{\U}{{\mathbf{U}}}
\newcommand{\V}{{\mathbf{V}}}
\newcommand{\X}{{\mathbf{X}}}
\newcommand{\Xp}{{\mathbf{X}^{(p)}}}
\newcommand{\Y}{{\mathbf{Y}}}
\newcommand{\Z}{{\mathbf{Z}}}

\newcommand{\KSch}{{K_{\mathrm{Sch}}}}
\newcommand{\Korth}{{K_{\mathrm{orth}}}}

\newcommand{\Beta}{{\mathrm{B}}}
\newcommand{\hM}{{\hat{\mathrm{M}}}}
\newcommand{\rM}{\mathrm{M}}
\newcommand{\Unif}{\operatorname{Unif}}

\newcommand{\tX}{{\tilde{\mathbf{X}}}}
\newcommand{\tY}{{\tilde{\mathbf{Y}}}}
\newcommand{\tBmu}{{\tilde{{\bolds{\mu}}}}}

\newcommand{\Bbeta}{{\bolds{\beta}}}
\newcommand{\Bmu}{{\bolds{\mu}}}

\newcommand{\hBb}{{\hat{\bolds{\beta}}}}
\newcommand{\hb}{{\hat{\beta}}}
\newcommand{\hY}{{\hat{\mathbf{Y}}}}
\newcommand{\hsig}{{\hat{\sigma}}}

\newcommand{\bPi}{\bolds{\Pi}}

\newcommand{\rank}{\operatorname{rank}}
\newcommand{\spanof}{\operatorname{span}}
\newcommand{\CI}{\operatorname{CI}}
\newcommand{\SPAR}{\mathrm{SPAR}}
\newcommand{\VIF}{\mathit{VIF}}

\newcommand{\argmin}{\mathop{\arg\min}}
\newcommand{\argmax}{\mathop{\arg\max}}

\makeatother

\begin{document}
\begin{frontmatter}

\title{Valid post-selection inference}
\runtitle{Valid post-selection inference}

\begin{aug}
\author[A]{\fnms{Richard} \snm{Berk}\ead[label=a1]{berk@wharton.upenn.edu}},
\author[A]{\fnms{Lawrence} \snm{Brown}\corref{}\thanksref{nsf}\ead[label=a2]{lbrown@wharton.upenn.edu}},
\author[A]{\fnms{Andreas} \snm{Buja}\thanksref{nsf}\ead[label=a3]{buja.at.wharton@gmail.com}},\\
\author[A]{\fnms{Kai} \snm{Zhang}\thanksref{nsf}\ead[label=a4]{zhangk@wharton.upenn.edu}}
\and
\author[A]{\fnms{Linda} \snm{Zhao}\thanksref{nsf}\ead[label=a5]{lzhao@wharton.upenn.edu}}
\runauthor{R. Berk et al.}
\affiliation{University of Pennsylvania}
\address[A]{Statistics Department\\
The Wharton School\\
University of Pennsylvania \\
471 Jon M. Huntsman Hall\\
Philadelphia, Pennsylvania 19104-6340 \\
USA\\
\printead{a1}\\
\hphantom{E-mail: }\printead*{a2}\\
\hphantom{E-mail: }\printead*{a3}\\
\hphantom{E-mail: }\printead*{a4}\\
\hphantom{E-mail: }\printead*{a5}} %adresu isvedimo komanda gale!
\end{aug}

\thankstext{nsf}{Supported in part by NSF Grant DMS-10-07657.}

% HISTORY:
\received{\smonth{7} \syear{2012}}
\revised{\smonth{12} \syear{2012}}

% ABSTRACT
%
\begin{abstract}
It is common practice in statistical data analysis to perform
data-driven variable selection and derive statistical inference from
the resulting model. Such inference enjoys none of the guarantees
that classical statistical theory provides for tests and confidence
intervals when the model has been chosen a priori. We propose to
produce valid ``post-selection inference'' by reducing the problem to
one of simultaneous inference and hence suitably widening conventional
confidence and retention intervals. Simultaneity is required for all
linear functions that arise as coefficient estimates in all submodels.
By purchasing ``simultaneity insurance'' for all possible submodels,
the resulting post-selection inference is rendered universally valid
under all possible model selection procedures. This inference is
therefore generally conservative for particular selection procedures,
but it is always less conservative than full Scheff\'{e} protection.
Importantly it does \textit{not} depend on the truth of the selected
submodel, and hence it produces valid inference even in wrong models.
We describe the structure of the simultaneous inference problem and
give some asymptotic results.
\end{abstract}

% KEYWORDS

\end{frontmatter}

%s1 #&#
\section{\texorpdfstring{Introduction: The problem with statistical inference after
model selection.}{Introduction: The problem with statistical inference after
model selection}}
\label{secintro}

Classical statistical theory grants validity of statistical tests and
confidence intervals assuming a wall of separation between the
selection of a model and the analysis of the data being modeled. In
practice, this separation rarely exists, and more often a model is
``found'' by a data-driven selection process. As a consequence
inferential guarantees derived from classical theory are invalidated.
Among model selection methods that are problematic for classical
inference, \textit{variable selection} stands out because it is regularly
taught, commonly practiced and highly researched as a technology.
Even though statisticians may have a general awareness that the
data-driven selection of variables (predictors, covariates) must
somehow affect subsequent classical inference from $F$- and $t$-based
tests and confidence intervals, the practice is so pervasive that it
appears in classical undergraduate textbooks on statistics such as
\citet{MooMcC03}.

The reason for the invalidation of classical inference guarantees is
that a data-driven variable selection process produces a model that is
itself stochastic, and this stochastic aspect is not accounted for by
classical theory. Models become stochastic when the stochastic
component of the data is involved in the selection process. (In
regression with fixed predictors the stochastic component is the
response.) Models are stochastic in a well-defined way when they are
the result of formal variable selection procedures such as stepwise or
stagewise forward selection or backward elimination or all-subset
searches driven by complexity penalties (such as $C_p$, AIC, BIC,
risk-inflation, LASSO$,\ldots$) or prediction criteria such as
cross-validation, or more recent proposals such as LARS and the Dantzig
selector; for an overview see, for example, \citet{HasTibFri09}.
Models are also stochastic but in an ill-defined way when they are
informally selected through visual inspection of residual plots or
normal quantile plots or other regression diagnostics. Finally, models
become stochastic in an opaque way when their selection is affected by
human intervention based on post-hoc considerations such as ``in
retrospect only one of these two variables should be in the model'' or
``it turns out the predictive benefit of this variable is too weak to
warrant the cost of collecting it.'' In practice, all three modes of
variable selection may be exercised in the same data analysis: multiple
runs of one or more formal search algorithms may be performed and
compared, the parameters of the algorithms may be subjected to
experimentation and the results may be critiqued with graphical
diagnostics; a~round of fine-tuning based on substantive deliberations
may finalize the analysis.

Posed so starkly, the problems with statistical inference after
variable selection may well seem insurmountable. At a minimum, one
would expect technical solutions to be possible only when a formal
selection algorithm is (1) well-specified (1a) in advance and
(1b) covering all eventualities, (2) strictly adhered to in the course
of data analysis and (3) not ``improved'' on by informal and post-hoc
elements. It may, however, be unrealistic to expect this level of
rigor in most data analysis contexts, with the exception of
well-conducted clinical trials. The real challenge is therefore to
devise statistical inference that is valid following \textit{any} type of
variable selection, be it formal, informal, post hoc or a combination
thereof. Meeting this challenge with a relatively simple proposal is
the goal of this article. This proposal for valid \textit{Po}st-\textit{S}election
\textit{I}nference, or ``\textit{PoSI}''
for short, consists of a large-scale family-wise error guarantee that
can be shown to account for all types of variable selection, including
those of the informal and post-hoc varieties. On the other hand, the
proposal\vadjust{\goodbreak} is no more conservative than necessary to account for
selection, and in particular it can be shown to be less conservative
than Scheff\'e's simultaneous inference.

The framework for our proposal is in outline as follows---details to
be elaborated in subsequent sections: we consider linear regression
with predictor variables whose values are considered fixed, and with a
response variable that has normal and homoscedastic errors. The
framework does not require that any of the eligible linear models is
correct, not even the full model, as long as a valid error estimate is
available. We assume that the selected model is the result of some
procedure that makes use of the response, but the procedure does not
need to be fully specified. A~crucial aspect of the framework
concerns the use and interpretation of the selected model: we assume
that, after variable selection is completed, the selected predictor
variables---and only they---will be relevant; all others will be
eliminated from further consideration. This assumption, seemingly
innocuous and natural, has critical consequences: it implies that
statistical inference will be sought for the coefficients of the
selected predictors only and in the context of the selected model
only. Thus the appropriate targets of inference are the best linear
coefficients within the selected model, where each coefficient is
adjusted for the presence of all other included predictors but not
those that were eliminated. Therefore the coefficient of an included
predictor generally requires inference that is specific to the model
in which it appears. Summarizing in a motto, a difference in
adjustment implies a difference in parameters and hence in inference.
The goal of the present proposal is therefore simultaneous inference
for all coefficients within all submodels. Such inference can be
shown to be valid following any variable selection procedure, be it
formal, informal, post hoc, fully or only partly specified.

Problems associated with post-selection inference were recognized long
ago, for example, by \citet{BueFed63}, \citet{Bro67},
\citet{Ols73}, \citet{Sen79}, \citet{SenSal87},
\citet{DijVel88}, \citet{Pot91}, \citet{Kab98}. More
recently specific problems have been the subject of incisive analyses
and criticisms by the ``Vienna School'' of P\"otscher, Leeb and
Schneider; see, for example, Leeb and P{\"o}tscher
(\citeyear{LeePot03,LeePot05,LeePot06N1,LeePot06N2,LeePot08N1,LeePot08N2,LeePot08N3}),
\citet{Pot06}, \citet{Lee06}, \citet{PotLee09},
P{\"o}tscher and Schneider (\citeyear{PotSch09,PotSch10,PotSch11}), as well as
\citet{KabLee06} and \citet{Kab09}. Important progress was
made by \citet{HjoCla03} and \citet{ClaHjo03}.

This article proceeds as follows: in
Section~\ref{secreinterpretation} we first develop the ``submodel
view'' of the targets of inference after model selection and contrast
it with the ``full model view'' (Section~\ref{secsubmodel}); we then
introduce assumptions with a view toward valid inference in ``wrong
models'' (Section~\ref{secassumptions}).
Section~\ref{secestimation} is about estimation and its targets from
the submodel point of view. Section~\ref{secvalidCIs} develops the
methodology for PoSI confidence intervals (CIs) and tests. After some
structural results for the PoSI problem in Section~\ref{secstruct},
we show in Section~\ref{secasymp}\vadjust{\goodbreak} that with increasing number of
predictors $p$ the width of PoSI CIs can range between the asymptotic
rates $O(\sqrt{\log{p}})$ and $O(\sqrt{p})$. We give examples for
both rates and, inspired by problems in sphere packing and covering,
we give upper bounds for the limiting constant in the $O(\sqrt{p})$
case. We conclude with a discussion in Section~\ref{secconc}. Some
proofs are deferred to the \hyperref[app]{Appendix}, and some elaborations to the
online Appendix in the supplementary material [\citet{Beretal}].

Computations will be described in a separate article.
Simulation-based methods yield satisfactory accuracy specific to a
design matrix up to $p \approx20$, while nonasymptotic universal
upper bounds can be computed for larger $p$.

%================================================================

%s2 #&#
\section{\texorpdfstring{Targets of inference and assumptions.}{Targets of inference and
assumptions}}
\label{secreinterpretation}

It is a natural intuition that model selection distorts inference by
distorting sampling distributions of parameter estimates: estimates in
selected models should tend to generate more type I errors than
conventional theory allows because the typical selection procedure
favors models with strong, hence highly significant predictors. This
intuition correctly points to a multiplicity problem that grows more
severe as the number of predictors subject to selection increases.
This is the problem we address in this article.

Model selection poses additional problems that are less obvious but no
less fundamental: there exists an ambiguity as to the role and meaning
of the parameters in submodels. In one view, the relevant parameters
are always those of the full model, hence the selection of a submodel
is interpreted as estimating the deselected parameters to be zero and
estimating the selected parameters under a zero constraint on the
deselected parameters. In another view, the submodel has its own
parameters, and the deselected parameters are not zero but
nonexistent. These distinctions are not academic as they imply
fundamentally different ideas regarding the targets of inference, the
measurement of statistical performance, and the problem of
post-selection inference. The two views derive from different
purposes of equations:
\begin{itemize}
\item Underlying the full model view of parameters is the use of a
full equation to describe a ``data generating'' mechanism for the
response; the equation hence has a causal interpretation.
\item Underlying the submodel view of parameters is the use of any
equation to merely describe association between predictor and
response variables; no data generating or causal claims are
implied.
\end{itemize}
In this article we address the latter use of equations. Issues
relating to the former use are discussed in the
online Appendix of the supplementary material [\citet{Beretal},
Section B.1].

%s2.1 #&#
\subsection{\texorpdfstring{The submodel interpretation of parameters.}{The submodel interpretation of
parameters}}
\label{secsubmodel}

In what follows we elaborate three points that set the submodel
interpretation of coefficients apart from the full model
interpretation, with important consequences for the rest of this
article:\vadjust{\goodbreak}
\begin{longlist}[(3)]
\item[(1)] The full model has no special status other than being
the repository of available predictors.
\item[(2)] The coefficients of excluded predictors are not zero;
they are not defined and therefore do not exist.
\item[(3)] The meaning of a predictor's coefficient depends on
which other predictors are included in the selected model.
\end{longlist}

(1) The full model available to the statistician often cannot be
argued to have special status because of inability to identify and
measure all relevant predictors. Additionally, even when a large and
potentially complete suite of predictors can be measured, there is
generally a question of predictor redundancy that may make it
desirable to omit some of the measurable predictors from the final
model. It is a common experience in the social sciences that models
proposed on theoretical grounds are found on empirical grounds to
have their predictors entangled by collinearities that permit little
meaningful statistical inference.
This situation is not limited to the social sciences: in gene
expression studies it may well occur that numerous sites have a
tendency to be expressed concurrently, hence as predictors in
disease studies they will be strongly confounded. The emphasis on
full models may be particularly strong in econometrics where there
is a ``notion that a longer regression \ldots has a causal
interpretation, while a shorter regression does not''
[\citet{AngPis09}, page 59]. Even in causal models,
however, there is a possibility that included adjuster variables
will ``adjust away'' some of the causal variables of interest.
Generally, in any creative observational study involving novel
predictors, it will be difficult a priori to exclude collinearities
that might force a rethinking of the predictors. In conclusion,
whenever predictor redundancy is a potential issue, it cannot a
priori be claimed that the full model provides the parameters of
primary interest.

(2) In the submodel interpretation of parameters, claiming that the
coefficients of deselected predictors are zero does not properly
describe the role of predictors. Deselected predictors have no role
in the submodel equation; they become no different than predictors
that had never been considered. The selected submodel becomes the
vehicle of substantive research irrespective of what the full model
was. As such the submodel stands on its own. This view is
especially appropriate if the statistician's task is to determine
which predictors are to be measured in the future.

(3) The submodel interpretation of parameters is deeply seated in how
we teach regression. We explain that the meaning of a regression
coefficient depends on which of the other predictors are included in
the model: ``the slope is the average difference in the response for
a unit difference in the predictor, \textit{at fixed levels of all
other predictors in the model}.'' This ``ceteris paribus'' clause
is essential to the meaning of a slope. That there is a difference
in meaning when there is a difference in covariates is most
drastically evident when there is a case of Simpson's paradox. For
example, if purchase likelihood of a high-tech gadget is predicted
from \textit{age}, it might be found against expectations that younger
people have lower purchase likelihood, whereas a regression on
\textit{age}
and \textit{income} might show that at fixed levels of income younger
people have indeed higher purchase likelihood. This case of
Simpson's paradox would be enabled by the expected positive
collinearity between \textit{age} and \textit{income}. Thus the marginal slope on
\textit{age} is distinct from the \textit{income}-adjusted slope on \textit{age} as the
two slopes answer different questions, apart from having opposite
signs. In summary, \textit{different models result in different
parameters with different meanings}.

Must we use the full model with both predictors? Not if \textit{income} data
is difficult to obtain or if it provides little improvement in $R^2$
beyond \textit{age}. The model based on \textit{age} alone cannot be said to be a
priori ``wrong.'' If, for example, the predictor and response
variables have jointly multivariate normal distributions, then every
linear submodel is ``correct.'' These considerations drive home, once
again, that sometimes no model has special status.

In summary, a range of applications call for a framework in which the
full model is not the sole provider of parameters, where rather each
submodel defines its own. The consequences of this view will be
developed in Section~\ref{secestimation}.

%s2.2 #&#
\subsection{\texorpdfstring{Assumptions, models as approximations, and error estimates.}{Assumptions, models as approximations, and error
estimates}}
\label{secassumptions}

We state assumptions for estimation and for the construction of valid
tests and CIs when fitting arbitrary linear equations. The main goal
is to prepare the ground for valid statistical inference after model
selection---\textit{not} assuming that selected models are correct.

We consider a quantitative response vector $\Y\in\reals^n$, assumed
random, and a full predictor matrix $\X=(\X_1,\X_2$,
$\ldots,\X_p)\in\reals^{n \times p}$, assumed fixed. We allow
$\X$ to be of nonfull rank, and $n$ and $p$ to be arbitrary. In
particular, we allow $n<p$. Throughout the article we let
%
%e2.1 #&#
\begin{equation}
\label{eqd} d \triangleq\rank(\X) = \dim\bigl(\spanof(\X)\bigr)\qquad
\mbox{hence } d
\le\min(n,p).
\end{equation}
Due to frequent reference we call $d=p$ ($\le n$) ``\textit{the classical
case}.''

It is common practice to assume the full model $\Y\sim
\mathcal{N}_n(\X\Bbeta,\sigma^2 \I)$ to be correct. In the present
framework, however, first-order correctness, $\E[\Y] = \X\Bbeta$,
will not be assumed. By implication, first-order correctness of any
submodel will not be assumed either. Effectively,
%
%e2.2 #&#
\begin{equation}
\label{eqmu} \Bmu\triangleq\E[\Y] \in\reals^n
\end{equation}
is allowed to be unconstrained and, in particular, need not reside in
the column space of $\X$. That is, the model given by $\X$ is allowed
to be ``first-order wrong,'' and hence we are, in a well-defined sense,
serious about G.E.P. Box's famous quote. What he calls ``wrong
models,'' we prefer to call ``approximations'': all predictor matrices
$\X$ provide approximations to $\Bmu$, some better than others, but
the degree of approximation plays no role in the clarification of
statistical inference.\vadjust{\goodbreak} The main reason for elaborating this point is
as follows: after model selection, the case for ``correct models'' is
clearly questionable, even for ``consistent model selection
procedures'' [\citet{LeePot03}, page 101]; but if
correctness of submodels is not assumed, it is only natural to abandon
this assumption for the full model also, in line with the idea that
the full model has no special status. As we proceed with estimation
and inference guarantees in the absence of first-order correctness we
will rely on assumptions as follows:
\begin{itemize}
\item For estimation (Section~\ref{secestimation}), we will
only need the existence of $\Bmu= \E[\Y]$.
\item For testing and CI guarantees (Section~\ref{secvalidCIs}),
we will make conventional second-order and distributional
assumptions,
%
%e2.3 #&#
\begin{equation}
\label{eqassumption} \Y\sim\cN\bigl(\Bmu, \sigma^2 \I\bigr).
\end{equation}
\end{itemize}
The assumptions (\ref{eqassumption}) of homoscedasticity and
normality are as questionable as first-order correctness, and we will
report elsewhere on approaches that avoid them. For now we follow the
vast model selection literature that relies on the technical
advantages of assuming homoscedastic and normal errors.

Accepting the assumption (\ref{eqassumption}), we address the issue
of estimating the error variance $\sigma^2$, because the valid tests
and CIs we construct require a valid estimate $\hsig^2$ of $\sigma^2$
that is independent of LS estimates. In the classical case, the most
common way to assert such an estimate is to assume that the full model
is first-order correct, $\Bmu= \X\Bbeta$ in addition to
(\ref{eqassumption}), in which case the mean squared residual (MSR)
$\hsig_F^2= \|\Y-\X\hBb\|^2 / (n-p)$ of the full model will do.
However, other possibilities for producing a valid estimate $\hsig^2$
exist, and they may allow relaxing the assumption of first-order
correctness:
\begin{itemize}
\item Exact replications of the response obtained under identical
conditions might be available in sufficient numbers. An estimate
$\hsig^2$ can be obtained as the MSR of the one-way ANOVA of the groups
of replicates.
\item In general,\vspace*{1pt} a larger linear model than the full model might be
considered as correct; hence $\hsig^2$ could be the MSR from this
larger model.
\item A different possibility is to use another dataset, similar to
the one currently being analyzed, to produce an independent estimate
$\hsig^2$ by whatever valid estimation method.
\item A special case of the preceding is a random split-sample
approach whereby one part of the data is reserved for producing
$\hsig^2$ and the other part for estimating coefficients, selecting
models and carrying out post-model selection inference.
\item A different type of estimate, $\hsig^2$, may be based on
considerations borrowed from nonparametric function estimation
[\citet{HalCar89}].
\end{itemize}
The purpose of pointing out these possibilities is to separate, at
least in principle, the issue of first-order model incorrectness from
the issue of error estimation under
assumption (\ref{eqassumption}). This separation puts the case
$n < p$ within our framework as the valid and independent estimation
of $\sigma^2$ is a problem faced by all ``$n < p$'' approaches.\vadjust{\goodbreak}

%================================================================

%s3 #&#
\section{\texorpdfstring{Estimation and its targets in submodels.}{Estimation and its targets in
submodels}}
\label{secestimation}

Following Section~\ref{secsubmodel}, the value and meaning of a
regression coefficient depends on what the other predictors in the
model are. An exception occurs, of course, when the predictors are
perfectly orthogonal, as in some designed experiments or in function
fitting with orthogonal basis functions. In this case a coefficient
has the same value and meaning across all submodels. This article is
hence a story of (partial) collinearity.

%s3.1 #&#
\subsection{\texorpdfstring{Multiplicity of regression coefficients.}{Multiplicity of regression
coefficients}}
\label{seccoeffs}

We will give meaning to LS estimators and their targets in the absence
of any assumptions other than the existence of $\Bmu= \E[\Y]$, which
in turn is permitted to be entirely unconstrained in $\reals^n$.
Besides resolving the issue of estimation in ``first-order wrong
models,'' the major purpose here is to elaborate the idea that the
slope of a predictor generates different parameters in different
submodels. As each predictor appears in $2^{p-1}$ submodels, the $p$
regression coefficients of the full model generally proliferate into a
plethora of as many as $p 2^{p-1}$ distinct regression coefficients
according to the submodels they appear in. To describe the situation
we start with notation.

To denote a submodel we use the (nonempty) index set $\rM=
\{j_1,j_2,\ldots,\break j_m\}\subset\rM_F = \{1,\ldots,p\}$ of the
predictors $\X_{j_i}$ in the submodel; the size of the submodel is
$m=|\rM|$ and that of the full model is $p=|\rM_F|$. Let $\X_\rM=
(\X_{j_1},\ldots,\X_{j_m})$ denote the $n \times m$ submatrix of $\X$
with columns indexed by~$\rM$. We will only allow submodels $\rM$ for
which $\X_\rM$ is of full rank,
\[
\rank(\X_\rM) = m \le d.
\]
We let $\hBb_\rM$ be the unique least squares estimate in $\rM$,
%
%e3.1 #&#
\begin{equation}
\label{eqestimate} \hBb_\rM= \bigl(\X_\rM^T
\X_\rM\bigr)^{-1}\X_\rM^T \Y.
\end{equation}

Now that $\hBb_\rM$ is an estimate, what is it estimating? Following
Section~\ref{secsubmodel}, we will not interpret $\hBb_\rM$ as
estimates of the full model coefficients and, more generally, of any
model other than $\rM$. Thus it is natural to ask that $\hBb_\rM$
define its own target through the requirement of unbiasedness,
%
%e3.2 #&#
\begin{equation}
\label{eqtarget} \Bbeta_\rM\triangleq\E[\hBb_\rM] =
\bigl(\X_\rM^T \X_\rM\bigr)^{-1}
\X_\rM^T \E[\Y] = \argmin_{\Bbeta' \in\reals^m} \bigl\| \Bmu-
\X_\rM\Bbeta' \bigr\| ^2.
\end{equation}
This definition requires no other assumption than the existence of
$\Bmu= \E[\Y]$. In particular there is no need to assume
first-order correctness of $\rM$ or $\rM_F$. Nor does it matter to what
degree $\rM$ provides a good approximation to $\Bmu$ in terms of
approximation error $\| \Bmu- \X_\rM\Bbeta_\rM\|^2$.

In the classical case $d = p \le n$, we can define the target of the
full-model estimate $\hBb= (\X^T \X)^{-1} \X^T \Y$ as a special case
of (\ref{eqtarget}) with $\rM=\rM_F$,
%
%e3.3 #&#
\begin{equation}
\label{eqtargetfull} \Bbeta\triangleq\E[\hBb] = \bigl(\X^T \X
\bigr)^{-1} \X^T \E[\Y].
\end{equation}
In the general (including the nonclassical) case, let $\Bbeta$ be any
(possibly nonunique) minimizer of $\|\Bmu- \X\Bbeta'\|^2$; the link\vadjust{\goodbreak}
between $\Bbeta$ and $\Bbeta_\rM$ is as follows:
%
%e3.4 #&#
\begin{equation}
\label{eqcontrast} \Bbeta_\rM= \bigl(\X_\rM^T
\X_\rM\bigr)^{-1}\X_\rM^T \X\Bbeta.
\end{equation}
Thus the target $\Bbeta_\rM$ is an estimable linear function of
$\Bbeta$, without first-order correctness assumptions. Equation
(\ref{eqcontrast}) follows from $\spanof(\X_\rM) \subset
\spanof(\X)$.

\textit{Notation}: to distinguish the regression coefficients of the
predictor $\X_j$ relative to the submodel it appears in, we write
$\beta_{j \cdot\rM} = \E[\hb_{j \cdot\rM}]$ for the components of
$\Bbeta_\rM= \E[\hBb_\rM]$ with $j \in\rM$. An important
convention is that indexes are always elements of the full model,
$j \in\{1,2,\ldots,p\} = \rM_F$, for what we call ``full model
indexing.''

%s3.2 #&#
\subsection{\texorpdfstring{Interpreting regression coefficients in first-order
incorrect models.}{Interpreting regression coefficients in first-order
incorrect models}}
\label{secinterp}

$\!\!\!$The regression coefficient $\beta_{j \cdot\rM}$ is conventionally
interpreted as the ``average difference in the response for a unit
difference in $X_j$, ceteris paribus in the model $\rM$.'' This
interpretation no longer holds when the assumption of first-order
correctness is given up. Instead, the phrase ``average difference in
the response'' should be replaced with the unwieldy phrase ``average
difference in the response approximated in the submodel $\rM$.'' The
reason is that the target of the fit $\hat{\Y}_\rM= \X_\rM\hBb
_\rM$
in the submodel $\rM$ is $\Bmu_\rM= \X_\rM\Bbeta_\rM$, hence in
$\rM$ we estimate unbiasedly not the true $\Bmu$ but its LS
approximation $\Bmu_\rM$.

A second interpretation of regression coefficients is in terms of
adjusted predictors: for $j \in\rM$ define the $\rM$-adjusted
predictor $\X_{j \cdot\rM}$ as the residual vector of the regression
of $\X_j$ on all other predictors in $\rM$. Multiple regression
coefficients, both estimates $\hb_{j \cdot\rM}$ and parameters
$\beta_{j \cdot\rM}$, can be expressed as simple regression
coefficients with regard to the $\rM$-adjusted predictors,
%
%e3.5 #&#
\begin{equation}
\label{eqadjusted} \hb_{j \cdot\rM} = \frac{\X_{j \cdot\rM}^T \Y}{\|\X
_{j \cdot
\rM}\|^2},\qquad \beta_{j \cdot\rM}
= \frac{\X_{j \cdot\rM}^T \Bmu}{\|\X
_{j \cdot\rM}\|^2}.
\end{equation}
The left-hand formula lends itself to an interpretation of $\hb_{j
\cdot\rM}$ in terms of the well-known leverage plot which shows
$\Y$ plotted against $\X_{j \cdot\rM}$ and the line with slope
$\hb_{j \cdot\rM}$. This plot is valid without first-order
correctness assumption.

A third interpretation can be derived from the second: to unclutter
notation let ${\mathbf{x}}= (x_i)_{i=1\ldots n}$ be any adjusted predictor
$\X_{j
\cdot\rM}$, so that $\hb= {\mathbf{x}}^T \Y/ \|{\mathbf{x}}\|^2$
and $\beta= {\mathbf{x}}^T
\Bmu/ \|{\mathbf{x}}\|^2$ are the corresponding $\hb_{j \cdot\rM}$ and
$\beta_{j \cdot\rM}$. Introduce (1) case-wise slopes through the
origin, both as estimates $\hb_{(i)} = Y_i/x_i$ and as parameters
$\beta_{(i)} = \mu_i/x_i$, and (2) case-wise weights $w_{(i)} = x_i^2
/ \sum_{i'=1\ldots n} x_{i'}^2$. Equations (\ref{eqadjusted}) are then
equivalent to the following:
\[
\hb= \sum_i w_{(i)}
\hb_{(i)},\qquad \beta= \sum_i
w_{(i)} \beta_{(i)}.
\]
Hence regression coefficients are weighted averages of case-wise
slopes, and this interpretation holds without first-order assumptions.\vadjust{\goodbreak}

%================================================================

%s4 #&#
\section{\texorpdfstring{Universally valid post-selection confidence intervals.}{Universally valid post-selection confidence
intervals}}
\label{secvalidCIs}

%s4.1 #&#
\subsection{\texorpdfstring{Test statistics with one error estimate for all submodels.}{Test statistics with one error estimate for all
submodels}}
\label{secerror}

We consider inference for $\hBb_\rM$ and its target $\Bbeta_\rM$.
Following Section~\ref{secassumptions} we require a normal
homoscedastic model for $\Y$, but we leave its mean $\Bmu= \E[\Y]$
entirely unspecified: $\Y\sim\cN(\Bmu, \sigma^2 \I)$. We then
have equivalently
\[
\hBb_\rM\sim\cN\bigl(\Bbeta_\rM, \sigma^2
\bigl(\X_\rM^T \X_\rM\bigr)^{-1}
\bigr) \quad\mbox{and}\quad \hb_{j \cdot\rM} \sim\cN\bigl(\beta_{j \cdot\rM},
\sigma^2 / \|\X_{j \cdot\rM}\|^2\bigr).
\]
Again\vspace*{1pt} following Section~\ref{secassumptions} we assume the
availability of a valid estimate $\hsig^2$ of $\sigma^2$ that is
independent of all estimates $\hb_{j \cdot\rM}$, and we further
assume $\hsig^2 \sim\sigma^2 \chi^2_{r} / r$ for $r$ degrees of
freedom. If the full model is assumed to be correct, $n>p$ and
$\hsig^2 = \hsig^2_F$, then $r = n - p$. In the limit $r
\rightarrow\infty$ we obtain $\hsig= \sigma$, the case of
known~$\sigma$, which will be used starting with
Section~\ref{secasymp}.

Let $t_{j \cdot\rM}$ denote a $t$-ratio for $\beta_{j \cdot\rM}$
that uses $\hsig$ irrespective of $\rM$,
%
%e4.1 #&#
\begin{equation}
\label{eqt} t_{j \cdot\rM} \triangleq\frac{\hb_{j \cdot\rM}-\beta_{j
\cdot\rM}} {
( (\X_{\rM}^T \X_{\rM})^{-1} )_{jj}^{{1/2}}
\hsig} = \frac{\hb_{j \cdot\rM}-\beta_{j \cdot\rM}} {
\hsig/ \|\X_{j \cdot\rM}\| } =
\frac{(\Y-\Bmu)^T \X_{j \cdot\rM}} {
\hsig\|\X_{j \cdot\rM}\| },
\end{equation}
where $(\cdots)_{jj}$ refers to the diagonal element corresponding to
$\X_j$. The quantity $t_{j \cdot\rM} = t_{j \cdot\rM}(\Y)$ has a
central $t$-distribution with $r$ degrees of freedom. Essential is
that the standard error estimate in the denominator of (\ref{eqt})
does \textit{not} involve the MSR $\hsig_\rM$ from the submodel $\rM$,
for two reasons:
\begin{itemize}
\item We do not assume that the submodel $\rM$ is first-order
correct; hence $\hsig_\rM^2$ would, in general, have a distribution
that is a multiple of a noncentral $\chi^2$ distribution with
unknown noncentrality parameter.
\item More disconcertingly, $\hsig_\rM^2$ would be the result of
selection, $\hsig_\hM^2$; see Section~\ref{secmodel-selection}.
Not much of real use is known about its distribution; see, for
example, \citet{Bro67} and \citet{Ols73}.
\end{itemize}
These problems are avoided by using one valid estimate $\hsig^2$ that
is independent of all submodels.

With this choice of $\hsig$, confidence intervals for $\beta_{j \cdot
\rM}$ take the form
%
%e4.2 #&#
\begin{eqnarray}
\label{eqjMCI} \CI_{j \cdot\rM}(K) &\triangleq& \bigl[ \hb_{j \cdot\rM}
\pm K
\bigl[ \bigl(\X_{\rM}^T \X_{\rM}
\bigr)^{-1} \bigr]_{jj}^{{1/2}} \hsig\bigr]
\nonumber\\[-8pt]\\[-8pt]
&=& \bigl[ \hb_{j \cdot\rM} \pm K \hsig/\|\X_{j \cdot\rM
}\| \bigr].
\nonumber
\end{eqnarray}
If $K = t_{r,1-\alpha/2}$ is the ${1 - \alpha/2}$ quantile of a
$t$-distribution with $r$ degrees of freedom, then the interval is
marginally valid with a $1 - \alpha$ coverage guarantee
\[
\P\bigl[\beta_{j \cdot\rM} \in\CI_{j \cdot\rM}(K)\bigr] \stackrel{(\ge)}
{=} 1-\alpha.
\]
This holds if the submodel $\rM$ is \textit{not} the result of variable
selection.

%s4.2 #&#
\subsection{\texorpdfstring{Model selection and its implications for parameters.}{Model selection and its implications for
parameters}}
\label{secmodel-selection}

In practice, the model $\rM$ tends to be the result of some form of
model selection that makes use of the stochastic component of the
data, which is the response vector~$\Y$ ($\X$ being\vspace*{1pt} fixed,
Section~\ref{secassumptions}). This model should therefore be
expressed as $\hM= \hM(\Y)$. In general we allow a variable
selection \textit{procedure} to be any (measurable) map
%
%e4.3 #&#
\begin{equation}
\label{eqmodelsel} \hM\dvtx  \Y\mapsto\hM(\Y),\qquad \reals^n \rightarrow\Mall,
\end{equation}
where $\Mall$ is the set of all full-rank submodels
%
%e4.4 #&#
\begin{equation}
\label{eqMall} \Mall\triangleq\bigl\{\rM| \rM\subset\{1,2,\ldots,p\},
\rank(
\X_\rM)=|\rM| \bigr\}.
\end{equation}
Thus the procedure $\hM$ is a discrete map that divides $\reals^n$
into as many as $|\cM_{\mathrm{all}}|$ different regions with shared outcome of
model selection.

Data dependence of the selected model $\hM$ has strong consequences:
\begin{itemize}
\item Most fundamentally, the selected model $\hM= \hM(\Y)$ is now
random. Whether the model has been selected by an algorithm or by
human choice, if the response $\Y$ has been involved in the
selection, the resulting model is a random object because it could
have been different for a different realization of the random
vector $\Y$.
\item Associated with the random model $\hM(\Y)$ is the parameter
vector of coefficients $\Bbeta_{\hM(\Y)}$, which is now randomly
chosen also:
\begin{itemize}
\item It has a random dimension $m(\Y) = |\hM(\Y)|$:
$ \Bbeta_{\hM(\Y)} \in\reals^{m(\Y)}$.
\item For any fixed\vspace*{1pt} $j$, it may or may not be the case
that $j \in\hM(\Y)$.
\item Conditional on $j \in\hM(\Y)$, the parameter
$\beta_{j \cdot\hM(\Y)}$ changes randomly as the adjuster
covariates in $\hM(\Y)$ vary randomly.
\end{itemize}
\end{itemize}
Thus the set of parameters for which inference is sought is random
also.

%s4.3 #&#
\subsection{\texorpdfstring{Post-selection coverage guarantees for confidence intervals.}{Post-selection coverage guarantees for confidence
intervals}}
\label{secpostselectionCIs}

With randomness of the selected model and its parameters in mind, what
is a desirable form of post-selection coverage guarantee for
confidence intervals? A natural requirement would be a $1 - \alpha$
confidence guarantee for the coefficients of the predictors that are
selected into the model,
%
%e4.5 #&#
\begin{equation}
\label{eqPoSIguarantee} \P\bigl[ \forall j \in\hM\dvtx  \beta_{j \cdot\hM}
\in
\CI_{j
\cdot\hM}(K) \bigr] \ge1-\alpha.
\end{equation}
Several points should be noted:
\begin{itemize}
\item The guarantee is family-wise for all selected predictors
$j \in\hM$, though the sense of ``family-wise'' is unusual
because $\hM=\hM(\Y)$ is random.
\item The guarantee has nothing to say about predictors
$j \notin\hM$ that have been deselected, regardless of the
substantive interest they might have. Predictors of overarching
interest should be protected from variable selection, and for
these one can use a modification of the PoSI approach which we
call ``PoSI1;'' see Section~\ref{secposi1}.\vadjust{\goodbreak}
\item Because predictor selection is random, $\hM= \hM(\Y)$, two
realized samples ${\mathbf{y}}^{(1)}, {\mathbf{y}}^{(2)} \in\reals
^n$ from $\Y$ may
result\vspace*{1pt} in different sets of selected predictors, $\hM({\mathbf{y}}^{(1)})
\neq\hM({\mathbf{y}}^{(2)})$. It would be a fundamental misunderstanding
to wonder whether the guarantee holds for both realizations.
Instead, the guarantee (\ref{eqPoSIguarantee}) is about the
\textit{procedure}
\[
\Y\mapsto\hsig(\Y),\qquad \hM(\Y), \hBb_{\hM(\Y)}(\Y) \mapsto
\CI_{j \cdot\hM}(K)\qquad (j \in\hM)
\]
for the long run of independent realizations of $\Y$ (by the LLN),
and not for any particular realizations ${\mathbf{y}}^{(1)}, {\mathbf
{y}}^{(2)}$. A
standard formulation used to navigate these complexities after a
realization ${\mathbf{y}}$ of $\Y$ has been analyzed is the
following: ``for
$j \in\hM$ we have $1 - \alpha$ \textit{confidence} that the
interval $\CI_{j \cdot\hM({\mathbf{y}})}(K)$ contains $\beta_{j
\cdot
\hM({\mathbf{y}})}$.''
\item Marginal\vspace*{-1pt} guarantees for individual predictors require some
care because $\beta_{j \cdot\hM}$ does not exist for
$j \notin\hM$. This makes\vspace*{-1pt} $\beta_{j \cdot\hM} \in\CI_{j
\cdot\hM}(K)$ an incoherent statement that does not define an
event. Guarantees are possible if the condition $j \in\hM$ is
added with a conjunction or is being conditioned on: The marginal
and conditional probabilities
\[
\P\bigl[ j \in\hM\mbox{ and } \beta_{j \cdot\hM} \in\CI_{j \cdot\hM}(K_{j \cdot})
\bigr] \quad\mbox{and}\quad \P\bigl[\beta_{j \cdot\hM} \in\CI_{j \cdot\hM}(K_{j
\cdot})|
j \in\hM\bigr],
\]
respectively, are both well-defined and can be the subject of coverage
guarantees; see the
online Appendix of the supplementary material [\citet{Beretal},
Section~B.4].
\end{itemize}

Finally, we note that the smallest constant $K$ that satisfies the
guarantee (\ref{eqPoSIguarantee}) is specific to the procedure $\hM$.
Thus different variable selection procedures would require different
constants. Finding procedure-specific constants is a challenge that
will be intentionally bypassed by the present proposals.

%s4.4 #&#
\subsection{\texorpdfstring{Universal validity for all selection procedures.}{Universal validity for all selection
procedures}}
\label{secuniv}

The ``PoSI'' procedure proposed here produces a constant $K$ that
provides universally valid post-selection inference \textit{for all model
selection procedures} $\hM$,
%
%e4.6 #&#
\begin{equation}
\label{eqbciuniv} \P\bigl[\beta_{j \cdot\hM} \in\CI_{j \cdot\hM}(K)
\ \forall
j \in\hM\bigr] \ge1-\alpha\qquad\forall\hM.
\end{equation}
Universal validity irrespective of the model selection procedure $\hM$
is a strong property that raises questions of whether the approach is
too conservative. There are, however, some arguments in its favor:

(1) Universal validity may be desirable or even essential for
applications in which the model selection procedure is not specified
in advance or for which the analysis involves some ad hoc elements
that cannot be accurately pre-specified. Even so, we should think of
the actually chosen model as part of a ``procedure'' $\Y\mapsto
\hM(\Y)$, and though the ad hoc steps are not specified for $\Y$ other
than the observed one, this is not a problem because our protection is
irrespective of what a specification might have been. This view also
allows data analysts to change their minds, to improvise and
informally decide in favor of a model other than that produced by a
formal selection procedure, or to experiment with multiple selection
procedures.

(2) There exists a model selection procedure that requires the full
strength of universally valid PoSI, and this procedure may not be
entirely unrealistic as an approximation to some types of data
analytic activities: ``significance hunting,'' that is, selecting that
model which contains the statistically most significant coefficient;
see Section~\ref{secspar}.

(3) There is a general question about the wisdom of proposing ever
tighter confidence and retention intervals for practical use when in
fact these intervals are valid only under tightly controlled
conditions. It might be realistic to suppose that much applied work
involves more data peeking than is reported in published articles.
With inference that is universally valid after any model selection
procedure, we have a way to establish which rejections are safe,
irrespective of unreported data peeking as part of selecting a model.

(4) Related to the previous point is the fact that today there is a
realization that a considerable fraction of published empirical work
is unreproducible or reports exaggerated effects; well known in this
regard is \citet{Ioa05}. A factor contributing to
this problem might well be liberal handling of variable selection and
absent accounting for it in subsequent inference.

%s4.5 #&#
\subsection{\texorpdfstring{Restricted model selection.}{Restricted model
selection}}
\label{secrestrictions}

The concerns over PoSI's conservative nature can be alleviated
somewhat by introducing a degree of flexibility to the PoSI problem
with regard to the universe of models being searched. Such
flexibility is additionally called for from a practical point of view
because it is not true that all submodels in $\Mall$ (\ref{eqMall})
are always being searched. Rather, the search is often limited in a
way that can be specified a priori, without involvement of $\Y$. For
example, a predictor of interest may be forced into the submodels of
interest, or there may be a restriction on the size of the submodels.
Indeed, if $p$ is large, a~restriction to a manageable set of
submodels is a computational necessity. In much of what follows we
can allow the universe $\cM$ of allowable submodels to be an (almost)
arbitrary but pre-specified nonempty subset of $\Mall$; w.l.o.g. we
can assume $\bigcup_{\rM\in\cM} \rM= \{1,2,\ldots,p\}$. Because we
allow only nonsingular submodels [see (\ref{eqMall})] we have
$|\rM| \le d$ $\forall\rM \in \cM$, where as always
$d = \rank(\X)$. Selection procedures are now maps
%
%e4.7 #&#
\begin{equation}
\label{eqmodelselrestr}
\hM\dvtx  \Y\mapsto\hM(\Y),\qquad \reals^n \rightarrow
\cM.
\end{equation}
The following are examples of model universes with practical relevance;
see also \citet{LeePot08N1}, Section 1.1,
Example 1.
\begin{longlist}[(5)]
\item[(1)] Submodels that contain the first $p'$ predictors ($1 \le
p' \le p$):
$\cM_1 = \{ \rM\in \Mall| \{1,2,\ldots,p'\}
\subset\rM\}$.
Classical: $|\cM_1| = 2^{p-p'}$. Example:
forcing an intercept into all models.\vadjust{\goodbreak}
\item[(2)] Submodels\vspace*{1pt} of size $m'$ or less (``sparsity option''):
$\cM_2 = \{ \rM\in\break \Mall| |\rM| \le m' \}$.
Classical: $|\cM_2| = {p\choose1} +\cdots+ {p\choose m'}$.\vspace*{1pt}
\item[(3)] Submodels with fewer than $m'$ predictors dropped from
the full model:
$\cM_3 = \{ \rM\in \Mall| |\rM| >
p-m' \}$. Classical: $|\cM_3| = |\cM_2|$.
\item[(4)] Nested models: $\cM_4 = \{ \{1,\ldots,j\} | j \in
\{1,\ldots,p\} \}$. $|\cM_4| = p$.
Example: selecting the degree up
to $p - 1$ in a polynomial regression.
\item[(5)] Models dictated by an ANOVA hierarchy of main effects and
interactions in a factorial design.
\end{longlist}
This list is just an indication of possibilities. In general, the
smaller the set $\tcM= \{(j,\rM) | j \in \rM\in \cM\}
$ is, the
less conservative the PoSI approach is, and the more computationally
manageable the problem becomes. With sufficiently strong
restrictions, in particular using the sparsity option (2) and assuming
the availability of an independent valid estimate $\hsig$, it is
possible to apply PoSI in certain nonclassical $p>n$ situations.

Further reduction of the PoSI problem is possible by pre-screening
adjusted predictors \textit{without the response} $\Y$. In a
fixed-design regression, any variable selection procedure that does
\textit{not} involve $\Y$ does \textit{not} invalidate statistical
inference. For example, one may decide not to seek inference for
predictors in submodels that impart a ``Variance Inflation Factor''
($\VIF$) above a user-chosen threshold: $\VIF_{j \cdot\rM} =
\|\X_j\|^2/\|\X_{j \cdot\rM}\|^2$ if $\X_j$ is centered, hence does
not make use of $\Y$, and elimination according to $\VIF_{j \cdot\rM}
> c$ does not invalidate inference.

%s4.6 #&#
\subsection{\texorpdfstring{Reduction of universally valid post-selection inference to
simultaneous inference.}{Reduction of universally valid post-selection inference to
simultaneous inference}}
\label{secreduction}

We show that universally valid post-selection inference
(\ref{eqbciuniv}) follows from simultaneous inference in the form of
family-wise error control for all parameters in all submodels. The
argument depends on the following lemma that may fall into the
category of the ``trivial but not immediately obvious.''

%le4.1 #&#
\begin{lem}[(``Significant triviality bound'')] \label{lemtrivial}
For any model
selection procedure $\hM\dvtx  \reals^n \rightarrow\cM$, the following
inequality holds for all $\Y\in\reals^n$:
\[
\max_{j \in\hM(\Y)} \bigl|t_{j \cdot\hM(\Y)}(\Y)\bigr| \le\max
_{\rM\in\cM} \max_{j \in\rM} \bigl|t_{j \cdot\rM}(\Y)\bigr|.
\]
\end{lem}

\begin{pf}
This is a special case of the triviality $f(\hM(\Y)) \le
\max_\rM f(\rM)$, where $f(\rM) = \max_{j \in\rM} |t_{j \cdot\rM
}(\Y)|$.
\end{pf}

The right-hand max-$|t|$ bound of the lemma is sharp in the sense that
there exists a variable selection procedure $\hM$ that attains the
bound; see Section~\ref{secspar}. Next we introduce the
$1-\alpha$ quantile of the right-hand max-$|t|$ statistic of the
lemma: let $K$ be the minimal value that satisfies
%
%e4.8 #&#
\begin{equation}
\label{eqKdef} \P\Bigl[\max_{\rM\in\cM} \max_{j \in\rM}
|t_{j \cdot\rM}| \le K \Bigr] \ge1-\alpha.\vadjust{\goodbreak}
\end{equation}
This value will be called ``the PoSI constant.'' It does not depend
on any model selection procedures, but it does depend on the design
matrix $\X$, the universe $\cM$ of models subject to selection, the
desired coverage $1-\alpha$, and the degrees of freedom $r$ in
$\hsig$, hence $K = K(\X,\cM,\alpha,r)$.

%th4.1 #&#
\begin{theorem} \label{thmsim}
For all model selection procedures $\hM\dvtx  \reals^n
\rightarrow\cM$ we have
%
%e4.9 #&#
\begin{equation}
\label{eqsim} \P\Bigl[ \max_{j \in\hM}|t_{j \cdot\hM}|\le K
\Bigr] \ge1-\alpha,
\end{equation}
where $K = K(\X,\cM,\alpha,r)$ is the PoSI constant.
\end{theorem}

This follows immediately from Lemma~\ref{lemtrivial}. Although
mathematically trivial we give the above the status of a theorem as it
is the central statement of the reduction of universal post-selection
inference to simultaneous inference. The following is just a
repackaging of Theorem~\ref{thmsim}:

%co4.1 #&#
\begin{cor} \label{corcov}
``Simultaneous post-selection confidence guarantees''
hold for any model selection procedure $ \hM\dvtx  \reals^n
\rightarrow \cM$,
%
%e4.10 #&#
\begin{equation}
\label{eqcov} \P\bigl[ \beta_{j \cdot\hM} \in\CI_{j \cdot\hM}(K)
\ \forall j
\in\hM\bigr] \ge1-\alpha,
\end{equation}
where $K = K(\X,\cM,\alpha,r)$ is the PoSI constant.
\end{cor}

Simultaneous inference provides strong family-wise error control,
which in turn translates to strong error control for tests following
model selection.

%co4.2 #&#
\begin{cor} \label{corstrong}
``Strong post-selection error control'' holds for any model selection
procedure $\hM\dvtx  \reals^n \rightarrow \cM$,
\[
\P\bigl[ \exists j \in\hM\dvtx  \beta_{j \cdot\hM} \neq0 \mbox{ and }
\bigl|t_{j \cdot\hM}^{(0)}\bigr|
> K \bigr] \le\alpha,
\]
where $K = K(\X,\cM,\alpha,r)$ is the PoSI constant and $t_{j
\cdot
\rM}^{(0)}$ is the $t$-statistic for the null hypothesis $\beta_{j
\cdot\rM} = 0$.
\end{cor}

The proof is standard; see the
online Appendix of the supplementary material [\citet{Beretal},
Section B.3].
The corollary states that, with probability $1-\alpha$, in a selected
model \textit{all} PoSI-significant rejections have detected true
alternatives.

%s4.7 #&#
\subsection{\texorpdfstring{Computation of the POSI constant.}{Computation of the POSI
constant}}

Several portions of the following treatment are devoted to a better
understanding of the structure and value of the POSI constant
$K(\X,\cM,\alpha,r)$. Except for very special choices it does not
seem possible to provide closed form expressions for its value.
However, the structural geometry and other properties to be described
later do enable a reasonably efficient computational algorithm.
R-code for computing the POSI\vadjust{\goodbreak} constant for small to moderate values of
$p$ is available on the authors' web pages. This code is accompanied
by a manuscript that will be published elsewhere describing the
computational algorithm and generalizations. For the basic setting
involving $\Mall$ the algorithm will conveniently provide values of
$K(\X,\Mall,\alpha,r)$ for matrices $\X$ of rank $\le20$, or slightly
larger depending on available computing speed and memory. It can also
be adapted to compute $K$ for some other families contained within
$\Mall$, such as some discussed in Section~\ref{secrestrictions}.

%s4.8 #&#
\subsection{\texorpdfstring{Scheff\'{e} protection.}{Scheff\'{e}
protection}}
\label{secscheffe}

Realizing the idea that the LS estimators in different submodels are
generally unbiased estimates of different parameters, we generated a
simultaneous inference problem involving up to $p 2^{ p-1}$ linear
contrasts $\beta_{j \cdot\rM}$. In view of the enormous number of
linear combinations for which simultaneous inference is sought, one
should wonder whether the problem is not best solved by Scheff\'{e}'s
method [\citet{Sch59}] which provides
simultaneous inference for \textit{all} linear combinations. To
accommodate rank-deficient $\X$, we cast Scheff\'e's result in terms
of $t$-statistics for arbitrary nonzero ${\mathbf{x}}\in\spanof(\X)$:
%
%e4.11 #&#
\begin{equation}
\label{eqt-general} t_{{\mathbf{x}}} \triangleq\frac{(\Y-\Bmu)^T
{\mathbf
{x}}}{\hsig\|{\mathbf{x}}\|}.
\end{equation}
The $t$-statistics in (\ref{eqt}) are obtained for ${\mathbf{x}}=
\X_{j \cdot\rM}$. Scheff\'e's guarantee is
%
%e4.12 #&#
\begin{equation}
\label{eqscheffe} \P\Bigl[\sup_{{\mathbf{x}}\in\spanof(\X)} |t_{\mathbf
{x}}| \le
\KSch\Bigr] = 1-\alpha,
\end{equation}
where the Scheff\'{e} constant is
%
%e4.13 #&#
\begin{equation}
\label{eqKSch} \KSch= \KSch(\alpha,d,r) = \sqrt{d \F_{d,r,1-\alpha}}.
\end{equation}
It provides an upper bound for \textit{all} PoSI constants:

%pr4.1 #&#
\begin{prop} \label{propscheffe}
$ K(\X,\cM,\alpha,r) \le \KSch(\alpha,d,r)
\ \forall\X, \cM, d = \rank(\X)$.
\end{prop}

Thus for $j \in\hM$ a parameter estimate $\hb_{j \cdot\hM}$ whose
$t$-ratio exceeds $\KSch$ in magnitude is universally safe from having
the rejection of ``$H_0\dvtx  \beta_{j \cdot\hM}=0$'' invalidated\vspace*{1pt} by
variable selection. The universality of the Scheff\'{e} constant is a
tip-off that it may be too loose for some predictor matrices $\X$, and
obtaining the sharper constant $K(\X)$ may be worthwhile. An
indication is given by the following comparison as $r \rightarrow
\infty$:
\begin{itemize}
\item For the Scheff\'{e} constant it holds $\KSch\sim\sqrt{d}$.
\item For orthogonal designs it holds $\Korth\sim\sqrt{2 \log{d}}$.
\end{itemize}
(For orthogonal designs see Section~\ref{secorthopt}.) Thus the PoSI
constant $\Korth$ is much smaller than $\KSch$. The large gap between
the two suggests that the Scheff\'{e} constant may be too conservative
at least in some cases. We will study certain nonorthogonal designs
for which the PoSI constant is $O(\sqrt{\log(d)})$ in\vadjust{\goodbreak}
Section~\ref{secexch}. On the other hand, the PoSI constant can
approach the order $O(\sqrt{d})$ of the Scheff\'{e} constant $\KSch$
as well, and we will study an example in Section~\ref{secrootp}.

Even though in this article we will give asymptotic results for $d=p
\rightarrow\infty$ and $r \rightarrow\infty$ only, we mention
another kind of asymptotics whereby $r$ is held constant while $d=p
\rightarrow\infty$: in this case $\KSch$ is in the order of the
product of $\sqrt{d}$ and the $1 - \alpha$ quantile of the
inverse-root-chi-square distribution with $r$ degrees of freedom. In
a similar way, the constant $\Korth$ for orthogonal designs is in the
order of the product of $\sqrt{2 \log{d}}$ and the $1 - \alpha$
quantile of the inverse-chi-square distribution with $r$ degrees of
freedom.

%s4.9 #&#
\subsection{\texorpdfstring{PoSI-sharp model selection---``SPAR.''}{PoSI-sharp model
selection---``SPAR''}}
\label{secspar}

There exists a model selection procedure that requires the full
protection of the simultaneous inference procedure (\ref{eqKdef}).
It is the ``significance hunting'' procedure that selects the model
containing the most significant ``effect'':
\[
\hM_{\SPAR}(\Y) \triangleq\argmax_{\rM\in\cM} \max
_{j \in\rM} \bigl| t_{j \cdot\rM}(\Y) \bigr|.
\]
We name this procedure ``SPAR'' for ``\textit{Single Predictor Adjusted
Regression}.'' It achieves equality with the ``significant
triviality bound'' in Lemma~\ref{lemtrivial} and is therefore the
worst case procedure for the PoSI problem. In the submodel
$\hM_{\SPAR}(\Y)$, the less significant predictors matter only in so
far as they boost the significance of the winning predictor by
adjusting it accordingly. This procedure ignores the quality of the
fit to $\Y$ provided by the model. While our present purpose is to
point out the existence of a selection procedure that requires full
PoSI protection, SPAR could be of practical interest when the analysis
is centered on strength of ``effects,'' not quality of model fit.

%s4.10 #&#
\subsection{\texorpdfstring{One primary predictor and controls---``PoSI\textit{1}.''}{One primary predictor and controls---``PoSI\textit{1}''}}
\label{secposi1}

Sometimes a regression analysis is centered on a predictor of
interest, $\X_j$, and on inference for its coefficient $\beta_{j
\cdot
\rM}$. The other predictors in $\rM$ act as controls, so their
purpose is to adjust the primary predictor for confounding effects and
possibly to boost the primary predictor's own ``effect.'' This
situation is characterized by two features:
\begin{itemize}
\item Variable selection is limited to models that contain the
primary predictor. We therefore define for any model universe
$\cM$ a sub-universe $\cM_{j \cdot}$ of models that contain the
primary predictor $\X_j$,
\[
\cM_{j \cdot} \triangleq\{ \rM| j \in\rM\in\cM\},
\]
so that for $\rM\in\cM$ we have $j \in\rM$ iff $\rM\in
\cM_{j
\cdot}$.
\item Inference is sought for the primary predictor $\X_j$ only,
hence the relevant test statistic is now $|t_{j \cdot\rM}|$ and
no longer $\max_{j \in\rM} |t_{j \cdot\rM}|$. The former
statistic is coherent because it is assumed that $j \in \rM$.
\end{itemize}
We call this the ``PoSI1'' situation in contrast to the unconstrained
PoSI situation. Similar to PoSI, PoSI1 starts with a ``significant
triviality bound'':\vadjust{\goodbreak}

%le4.2 #&#
\begin{lem}[(``Primary predictor's significant triviality bound'')] \label{lemtrivial1}
For a fixed predictor $\X_j$ and model selection procedure $\hM\dvtx
\reals^n \rightarrow\cM_{j \cdot}$, it holds that
\[
\bigl|t_{j \cdot\hM(\Y)}(\Y)\bigr| \le\max_{\rM\in\cM_{j \cdot}} \bigl|t_{j \cdot\rM
}(\Y)\bigr|.
\]
\end{lem}

For a ``proof,'' the only thing to note is $j \in \hM(\Y)$ by the
assumption $\hM(\Y) \in \cM_{j \cdot}$. We next define the
``PoSI1'' constant for the predictor $\X_j$ as the $1 - \alpha$
quantile of the max-$|t|$ statistic on the right-hand side of the
lemma: let
$K_{j \cdot} = K_{j \cdot}(\X,\cM,\alpha,r)$ be the minimal value that
satisfies
%
%e4.14 #&#
\begin{equation}
\label{eqKjdef} \P\Bigl[ \max_{\rM\in\cM_{j \cdot}} |t_{j \cdot\rM}|
\le
K_{j
\cdot} \Bigr] \ge1 - \alpha.
\end{equation}
Importantly, this constant is dominated by the general PoSI constant,
%
%e4.15 #&#
\begin{equation}
\label{eqKvsKj} K_{j \cdot}(\X,\cM,\alpha,r) \le K(\X,\cM,\alpha,r)
\end{equation}
for the obvious reason that the present max-$|t|$ is smaller than the
general PoSI max-$|t|$ due to $\cM_{j \cdot} \subset\cM$ and the
restriction of inference to $\X_j$. The constant $K_{j \cdot}$
provides the following ``PoSI1'' guarantee shown as the analog of
Theorem~\ref{thmsim} and Corollary~\ref{corcov} folded into one:

%th4.2 #&#
\begin{theorem} \label{thmsim1r}
Let $\hM\dvtx  \reals^n \rightarrow\cM_{j \cdot}$ be a
selection procedure that always includes the predictor $\X_j$ in the
model. Then we have
%
%e4.16 #&#
\begin{equation}
\label{eqsimjr} \P\bigl[ |t_{j \cdot\hM}| \le K_{j \cdot} \bigr] \ge1 - \alpha,
\end{equation}
and accordingly we have the following post-selection confidence
guarantee:
%
%e4.17 #&#
\begin{equation}
\label{eqcovjr} \P\bigl[ \beta_{j \cdot\hM} \in\CI_{j \cdot\hM}(K_{j
\cdot
})
\bigr] \ge1-\alpha.
\end{equation}
\end{theorem}

Inequality (\ref{eqsimjr}) is immediate from
Lemma~\ref{lemtrivial1}. The ``triviality bound'' of the lemma is
attained by the following variable selection procedure which we name
``SPAR1'':
%
%e4.18 #&#
\begin{equation}
\label{eqspar1} \hM_{j \cdot}(\Y) \triangleq\argmax_{\rM\in\cM_{j \cdot
}} \bigl|
t_{j \cdot\rM}(\Y) \bigr|.
\end{equation}
It is a potentially realistic description of some data analyses when a
predictor of interest is determined a priori, and the goal is to
optimize \textit{this} predictor's ``effect.'' This procedure requires
the full protection of the PoSI1 constant~$K_{j \cdot}$.

In addition to its methodological interest, the PoSI1 situation
addressed by Theorem~\ref{thmsim1r} is of theoretical interest: even
though the PoSI1 constant~$K_{j \cdot}$ is dominated by the
unrestricted PoSI constant $K$, we will construct in
Section~\ref{secrootp} an example of predictor matrices for which the
PoSI1 constant increases at the Scheff\'{e} rate and is asymptotically
more than 63\% of the Scheff\'e constant $\KSch$. It follows that
near-Scheff\'e protection can be needed even for SPAR1 variable
selection.

%================================================================

%s5 #&#
\section{\texorpdfstring{The structure of the PoSI problem.}{The structure of the PoSI
problem}}
\label{secstruct}

%s5.1 #&#
\subsection{\texorpdfstring{Canonical coordinates.}{Canonical
coordinates}}
\label{seccancoord}

We can reduce the dimensionality of the PoSI problem from $n \times p$
to $d \times p$, where $d = \rank(X) \le\min(n,p)$, by introducing
Scheff\'e's canonical coordinates. This reduction is important both
geometrically and computationally because the PoSI coverage problem
really takes place in the column space of $\X$.

\begin{definition*}
Let $\Q= ({\mathbf{q}}_1,\ldots,{\mathbf {q}}_d) \in\reals^{n \times
d}$ be any orthonormal basis of the column space of $\X$. Note that
$\hY= \Q\Q^T \Y$ is the orthogonal projection of $\Y$ onto the column
space of $\X$ even if $\X$ is not of full rank. We call $\tX= \Q^T
\X\in\reals^{d \times p}$ and $\tY= \Q^T \hY\in \reals^d$ canonical
coordinates of $\X$ and $\hY$.
\end{definition*}

We extend the notation $\X_\rM$ for extraction of subsets of columns
to canonical coordinates $\tX_\rM$. Accordingly slopes obtained from
canonical coordinates will be denoted by $\hBb_\rM(\tX,\tY) =
(\tX_\rM^T \tX_\rM)^{-1} \tX_\rM^T \tY$ to distinguish them from the
slopes obtained from the original data $\hBb_\rM(\X,\Y) = (\X_\rM^T
\X_\rM)^{-1} \X_\rM^T \Y$, if only to state in the following
proposition that they are identical.

%pr5.1 #&#
\begin{prop}\label{propcancoord}
Properties of canonical coordinates:
\begin{longlist}[(7)]
\item[(1)] $\tY= \Q^T \Y$.
\item[(2)] $\tX_\rM^T \tX_\rM= \X_\rM^T \X_\rM$ and $\tX_\rM
^T \tY= \X_\rM^T \Y$.
\item[(3)] $\hBb_\rM(\tX,\tY) = \hBb_\rM(\X,\Y)$ for all
submodels $M$.
\item[(4)] $\tY\sim\cN(\tBmu, \sigma^2 \I_d)$, where $\tBmu= \Q
^T \Bmu$.
\item[(5)] $\tX_{j \cdot\rM} = \Q^T \X_{j \cdot\rM}$, where $j
\in\rM$
and $\tX_{j \cdot\rM} \in\reals^d$ is the residual vector of the
regression of $\tX_j$ onto the other columns of $\tX_\rM$.
\item[(6)] $t_{j \cdot M} =
(\hb_{j \cdot M}(\tX,\tY) - \beta_{j \cdot M})/(\hsig/\|\tX_{j
\cdot\rM}\|)$.
\item[(7)] In the classical case $d = p$, $\tX$ can be chosen to be
an upper triangular or a symmetric matrix.
\end{longlist}
\end{prop}

The proofs of (1)--(6) are elementary. As for
(7), an upper triangular $\tX$ can be obtained from a
QR-decomposition based on a Gram--Schmidt procedure, $\X= \Q\R$,
$\tX
= \R$. A symmetric $\tX$ is obtained from a singular value
decomposition, $\X= \U\D\V^T$, $\Q= \U\V^T$, $\tX= \V\D\V^T$.

Canonical coordinates allow us to analyze the PoSI coverage problem in~$\reals^d$. In what follows we will freely assume that all objects
are rendered in canonical coordinates and write $\X$ and $\Y$ for
$\tX$ and $\tY$, implying that the predictor matrix is of size $d
\times p$ and the response is of size $d \times1$.

%s5.2 #&#
\subsection{\texorpdfstring{PoSI coefficient vectors in canonical coordinates.}{PoSI coefficient vectors in canonical
coordinates}}
\label{secPoSIcanonical}

We simplify the PoSI coverage problem (\ref{eqKdef}) as follows: due
to pivotality of $t$-statistics, the problem is invariant under
translation of $\Bmu$ and rescaling of $\sigma$; see equation
(\ref{eqt}). Hence it suffices to solve coverage problems for\vadjust{\goodbreak}
$\Bmu=\0$ and $\sigma=1$. In canonical coordinates this implies
$\E[\tY]=\0_d$, hence $\tY\sim\cN(\0_d,\I_d)$. For this reason we
use the more familiar notation $\Z$ instead of $\tY$. The random
vector $\Z/\hsig$ has a $d$-dimensional $t$-distribution with $r$
degrees of freedom, and any linear combination ${\mathbf{u}}^T \Z
/\hsig$ with a
unit vector ${\mathbf{u}}$ has a one-dimensional $t$-distribution.
Letting $\X_{j
\cdot\rM}$ be the adjusted predictors in canonical coordinates,
estimates (\ref{eqadjusted}) and their $t$-statistics (\ref{eqt})
simplify to
%
%e5.1 #&#
\begin{equation}
\label{eqz1} \hb_{j \cdot\rM} = \frac{\X_{j \cdot\rM}^T \Z}{\|\X_{j
\cdot
\rM}\|^2} = \l_{j \cdot\rM}^T
\Z,\qquad t_{j \cdot\rM} = \frac{\X_{j \cdot\rM}^T \Z}{\|\X_{j \cdot
\rM}\| \hsig} = \ls_{j \cdot\rM}^T \Z/
\hsig,
\end{equation}
which are linear functions of $\Z$ and $\Z/\hsig$, respectively,
with ``PoSI coefficient vectors'' $\l_{j \cdot\rM}$ and $\ls_{j
\cdot\rM}$ that equal $\X_{j \cdot\rM}$ up to scale
%
%e5.2 #&#
\begin{equation}
\label{eqposivecs} \l_{j \cdot\rM} \triangleq\frac{\X_{j \cdot\rM}}{\|
\X_{j
\cdot\rM}\|^2},\qquad
\ls_{j \cdot\rM} \triangleq\frac{\l_{j \cdot\rM}}{\|\l_{j
\cdot\rM}\|} = \frac{\X_{j \cdot\rM}}{\|\X_{j \cdot\rM}\|}.
\end{equation}
As we now operate in canonical coordinates, we have $\l_{j \cdot\rM}
\in \reals^d$ and $\ls_{j \cdot\rM} \in S^{d-1}$, the unit
sphere in $\reals^d$. To complete the structural description of the
PoSI problem we let
%
%e5.3 #&#
\begin{equation}
\label{eqL} \cL(\X,\cM) \triangleq\{ \ls_{j \cdot\rM} | j \in\rM\in
\cM\}
\subset S^{d-1}.
\end{equation}
If $\cM= \Mall$, we omit the second argument and write $\cL(\X)$.

%pr5.2 #&#
\begin{prop} \label{propmaxlin}
The PoSI problem (\ref{eqKdef}) is equivalent to a $d$-dimen\-sional
coverage problem for linear functions of the multivariate
$t$-vector $\Z/\hsig$,
%
%e5.4 #&#
\begin{equation}
\label{eqtcoverage} \P\Bigl[\max_{\rM\in\cM}\max
_{j \in\rM} |t_{j \cdot\rM}| \leq K \Bigr] = \P\Bigl[\max
_{\ls\in\cL(\X,\cM)} \bigl|\ls^T\Z/\hsig\bigr| \leq K \Bigr] \stackrel{(
\ge)} {=} 1 - \alpha.
\end{equation}
\end{prop}

%s5.3 #&#
\subsection{\texorpdfstring{Orthogonalities of PoSI coefficient vectors.}{Orthogonalities of PoSI coefficient
vectors}}
\label{secorthgeo}

The set $\cL(\X,\cM)$ of unit vectors $\ls_{j \cdot\rM}$ has
interesting geometric structure which is the subject of this and the
next subsection. The following proposition (proof in
Appendix~\ref{appproporth}) elaborates the fact that $\ls_{j \cdot
\rM}$ is essentially the predictor vector $\X_j$ orthogonalized with
regard to the other predictors in the model $\rM$. Vectors will
always be assumed in canonical coordinates and hence $d$-dimensional.

%pr5.3 #&#
\begin{prop} \label{proporth}
Orthogonalities in $\cL(\X,\cM)$: the following statements hold
assuming that the models referred to are in $\cM$ (hence are of full
rank).
\begin{longlist}[(4)]
\item[(1)] Adjustment properties:
\[
\ls_{j \cdot\rM} \in\spanof\{ \X_j | j \in\rM\}\quad
\mbox{and}\quad
\ls_{j \cdot\rM} \perp\X_{j'} \qquad \mbox{for $j \neq j'$ both $\in\rM$}.
\]
\item[(2)] The following vectors form an orthonormal ``Gram--Schmidt''
series:
\[
\{\ls_{1 \cdot\{1\}}, \ls_{2 \cdot\{1,2\}}, \ls_{3 \cdot
\{1,2,3\}},\ldots,
\ls_{d \cdot\{1,2,\ldots,d\}} \}.
\]
Other series are obtained using $(j_1,j_2,\ldots,j_d)$ in place of
$(1,2,\ldots,d)$.
\item[(3)] Vectors $\ls_{j \cdot\rM}$ and $\ls_{j' \cdot\rM'}$ are
orthogonal if $\rM\subset\rM'$, $j \in\rM$ and \mbox{$j' \in
\rM'
\setminus\rM$}.
\item[(4)] Classical case $d = p$ and $\cM= \Mall$: each vector
$\ls_{j \cdot\rM}$ is orthogonal to $(p - 1) 2^{p-2}$ vectors
$\ls_{j' \cdot\rM'}$ (not all of which may be distinct).\vadjust{\goodbreak}
\end{longlist}
\end{prop}

The cardinality of orthogonalities in the classical case and
$\cM= \Mall$ is as follows: if the predictor vectors $\X_j$ have no
orthogonal pairs among them, then $|\cL(\X)| = p 2^{p-1}$. If
there exist orthogonal pairs, then $|\cL(\X)|$ is less. For example,
if there exists exactly one orthogonal pair, then $|\cL(\X)| =
(p - 1) 2^{p-1}$. When $\X$ is a fully orthogonal design, then
$|\cL(\X)| = p$.

%s5.4 #&#
\subsection{\texorpdfstring{The PoSI polytope.}{The PoSI polytope}}
\label{secpolytope}

Coverage problems can be framed geometrically in terms of probability
coverage of polytopes in $\reals^d$. For the PoSI problem the
polytope with half-width $K$ is defined by
%
%e5.5 #&#
\begin{equation}
\label{eqPoSIpolytope} \bPi_K = \bPi_K(\X,\cM)
\triangleq\bigl\{ {\mathbf{z}}\in\bbR^d | \bigl|\ls^T {
\mathbf{z}}\bigr| \le K, \forall\ls\in\cL(\X,\cM) \bigr\},
\end{equation}
henceforth called the ``PoSI polytope.'' The PoSI coverage problem
(\ref{eqtcoverage}) is equivalent to calibrating $K$ such that
\[
\P[\Z/\hsig\in\bPi_K] = 1-\alpha.
\]
The simplest case of a PoSI polytope, for $d = p = 2$, is
illustrated in Figure 1 %%%%%%%%%%\ref{figpoly2}
in the
online Appendix of the supplementary material [\citet{Beretal},
Section~B.7]. More general polytopes\vspace*{1pt} are
obtained for arbitrary sets $\cL$ of unit vectors, that is, subsets
$\cL\subset S^{d-1}$ of the unit sphere in $\reals^d$. For the
special case $\cL= S^{d-1}$ the ``polytope'' is the ``Scheff\'{e}
ball'' with coverage $\sqrt{d\F_{d,r}} \rightarrow\sqrt{\chi^2_d}$ as
$r \rightarrow\infty$:
\[
\B_K \triangleq\bigl\{ {\mathbf{z}}\in\bbR^d | \|{
\mathbf{z}}\| \le K \bigr\},\qquad \P[ \Z/\hsig\in\B_K ] = F_{\F_{d,r}}
\bigl(K^2/d\bigr).
\]

Many properties of the polytopes $\bPi_K$ are not specific to PoSI
because they hold for polytopes (\ref{eqPoSIpolytope}) generated by
simultaneous inference problems for linear functions with arbitrary
sets $\cL$ of unit vectors. These polytopes:
\begin{longlist}[(5)]
\item[(1)]form scale families of geometrically similar bodies:
$\bPi_K = K\bPi_1$;
\item[(2)] are point symmetric about the origin: $\bPi_K =
-\bPi_K$;
\item[(3)] contain the Scheff\'{e} ball: $\B_K \subset\bPi_K$;
\item[(4)] are intersections of ``slabs'' of width $2K$:
\[
\bPi_K = \bigcap_{\ls\in\cL} \bigl\{{
\mathbf{z}}\in\bbR^d | \bigl|{\mathbf{z}}^T\ls\bigr| \le K \bigr\};
\]
\item[(5)] have $2 |\cL|$ faces (assuming $\cL\cap-\cL= \varnothing$),
and each face is tangent to the Scheff\'{e} ball $\B_K$ with
tangency points $\pm K \ls$ $(\ls\in\cL)$.
\end{longlist}
Specific to PoSI are the orthogonalities described in
Proposition~\ref{proporth}.

%s5.5 #&#
\subsection{\texorpdfstring{PoSI optimality of orthogonal designs.}{PoSI optimality of orthogonal
designs}}
\label{secorthopt}

In orthogonal designs, adjustment has no effect: $\X_{j \cdot
\rM} = \X_j$ for all $j \in \rM$, hence $\ls_{j \cdot
\rM} = \X_j/\|\X_j\|$ and $\cL(\X,\cM) = \{ \X_1/\|\X_1\|,\ldots,
\X_p/\|\X_p\| \}$. The polytope $\bPi_K$ is therefore a hypercube.
This observation implies an optimality property of orthogonal designs
if the submodel universes $\cM$ are sufficiently rich to force
$\cL(\X,\cM)$ to contain an orthonormal basis\vadjust{\goodbreak} of $\reals^d$: the
polytope generated by an orthonormal basis is a hypercube; hence the
polytope $\bPi_K(\X,\cM)$ is contained in this hypercube; thus
$\bPi_K(\X,\cM)$ has maximal extent if and only if it is equal to
this hypercube,
which is the case if and only if $\cL(\X,\cM)$ is this orthonormal
basis and
nothing more; that is, $\X$ is an orthogonal design. A~simple
sufficient condition for $\cM$ to grant the existence of an
orthonormal basis in $\cL(\X,\cM)$ is the existence of a maximal
nested sequence of submodels such as $\{1\}$,
$\{1,2\},\ldots,\{1,2,\ldots,d\}$ in $\cM$. It follows according to
item (2) in Proposition~\ref{proporth} that there exists an
orthonormal Gram--Schmidt basis in $\cL(\X,\cM)$. We summarize:

%pr5.4 #&#
\begin{prop} \label{proporthopt}
Among predictor matrices with $\rank(\X) = d$ and mod\-el
universes $\cM$ that contain at least one maximal nested sequence of
submodels, orthogonal designs with $p = d$ columns yield:
\begin{itemize}
\item the maximal coverage probability $\P[\Z/\hsig\in\bPi_K]$ for
fixed $K$ and
\item the minimal PoSI constant $K$ satisfying $\P[\Z/\hsig\in\bPi
_{K}] =
1-\alpha$ for fixed~$\alpha$, $\inf_{\rank(X)=d} K(\X,\cM,\alpha,r) =
\Korth(\alpha,d,r)$.
\end{itemize}
\end{prop}

The proposition holds not only for multivariate $t$-vectors and their
Gaussian limits but for arbitrary spherically symmetric distributions.
Optimality of orthogonal designs translates to optimal asymptotic
behavior of their constant $K(\X,\alpha)$ for large $d$:

%pr5.5 #&#
\begin{prop} \label{proporthasy}
Consider the Gaussian limit $r \rightarrow\infty$. For $\X$ and
$\cM$ as in Proposition~\ref{proporthopt}, the asymptotic lower bound
for the constant $K$ as $d \rightarrow\infty$ is attained for
orthogonal designs for which the asymptotic rate is
\[
\inf_{\rank(\X)=d} K(\X,\cM,\alpha) = \Korth(d,\alpha) = \sqrt{2
\log{d}}+o(d).
\]
\end{prop}

By Proposition~\ref{proporthopt} the PoSI problem is bounded below by
orthogonal designs, and by Proposition~\ref{propscheffe} it is
loosely bounded above by the Scheff\'{e} ball (both for all $\alpha$,
$d$, and $r$). The question of how close to the Scheff\'{e} bound
PoSI problems can get for $r \rightarrow\infty$ will occupy us in
Section~\ref{secrootp}. Unlike the infimum problem, the supremum
problem does not appear to have a unique optimizing design $\X$
uniformly in $\alpha$, $d$ and $r$.

%s5.6 #&#
\subsection{\texorpdfstring{A duality property of PoSI vectors.}{A duality property of PoSI
vectors}}
\label{secdual}

In the classical case $d = p$ and $\cM= \Mall$ there exists a
duality for PoSI vectors $\cL(\X)$ which we will use in
Section~\ref{secexch} below but which is also of independent
interest. Some preliminaries: letting $\rM_F = \{1,2,\ldots,p\}$ be the
full model, we observe that the (unnormalized) PoSI vectors $\l_{j
\cdot\rM_F} = \X_{j \cdot\rM_F} / \|\X_{j \cdot\rM_F}\|^2$ form
the rows of the matrix $(\X^T \X)^{-1} \X^T$; see (\ref{eqadjusted})
and (\ref{eqcontrast}). In a change of perspective, we interpret
the transpose matrix
\[
\X^* = \X\bigl(\X^T \X\bigr)^{-1}
\]
as a predictor matrix, to be called the ``dual design'' of $\X$. It
is also of size $p \times p$ in canonical coordinates, and its columns
are the PoSI vectors $\l_{j \cdot\rM_F}$. It turns out that $\X^*$
and $\X$ pose identical PoSI problems if $\cM= \Mall$:

%th5.1 #&#
\begin{theorem} \label{thmdual}
$\cL(\X^*) = \cL(\X), \bPi_K(\X^*) = \bPi_K(\X), K(\X^*)
= K(\X) $.
\end{theorem}

Recall that $\cL(\X)$ and $\cL(\X^*)$ contain the normalized versions
of the respective adjusted predictor vectors. The theorem follows
from the following lemma which establishes the identities of vectors
between $\cL(\X^*)$ and $\cL(\X)$. We extend obvious notation from
$\X$ to $\X^*$ as follows:
\[
\X_j^* = \l^*_{j \cdot\{j\}} = \l_{j \cdot\rM_F}.
\]
Submodels for $\X^*$ will be denoted $\rM^*$, but they, too, will be
given as subsets of $\{1,2,\ldots,p\}$ which, however, refer to columns
of $\X^*$. Finally, the normalized version of $\l^*_{j \cdot\rM^*}$
will be written as $\ls^*_{j \cdot\rM^*}$.

%le5.1 #&#
\begin{lem} \label{lemdual}
For two submodels $\rM$ and $\rM^*$ that satisfy $\rM\cap\rM^* = \{
j\}$
and $\rM\cup\rM^* = \rM_F$, we have
\[
\ls^*_{j \cdot\rM^*} = \ls_{j \cdot\rM},\qquad \bigl\|\l^*_{j \cdot\rM^*}\bigr\| \|
\l_{j \cdot\rM}\| = 1.
\]
\end{lem}

The proof is in Appendix~\ref{applemdual}. The assertion about
norms is really only needed to exclude collapse of $\l^*_{j \cdot
\rM^*}$ to zero.

A special case arises when the predictor matrix (in canonical
coordinates) is chosen to be symmetric according to
Proposition~\ref{propcancoord}(7): if $\X^T = \X$, then $\X^*
= \X(\X^T\X)^{-1}=\X^{-1}$, and hence:

%co5.1 #&#
\begin{cor} \label{symdualK}
If $\X$ is symmetric in canonical coordinates, then
\[
\cL\bigl(\X^{-1}\bigr) = \cL(\X),\qquad \bPi_K\bigl(
\X^{-1}\bigr) = \bPi_K(\X) \quad\mbox{and}\quad K\bigl(\X^{-1}
\bigr) = K(\X).
\]
\end{cor}

%================================================================

%s6 #&#
\section{\texorpdfstring{Illustrative examples and asymptotic results.}{Illustrative examples and asymptotic
results}}
\label{secasymp}

We consider examples in the classical case $d = p$ and
$\cM= \Mall$. Also, we work with the Gaussian limit $r \rightarrow
\infty$, that is, $\sigma^2$ known, and w.l.o.g. $\sigma^2=1$.

%s6.1 #&#
\subsection{\texorpdfstring{Example 1: Exchangeable designs.}{Example 1: Exchangeable
designs}}
\label{secexch}

In exchangeable designs all pairs of predictor vectors enclose
the same angle. In canonical coordinates a convenient parametrization
of a family of symmetric exchangeable designs is
%
%e6.1 #&#
\begin{equation}
\label{eqex} \Xp(a) = \I_p + a \E_{p \times p},
\end{equation}
where $ -1/p < a < \infty$, and $\E_{p \times p}$ is a matrix with all
entries equal to $1$. The range restriction on $a$ assures that $\Xp$
is positive definite. We will write $\X= \Xp= \X(a)=\Xp(a)$
depending on which parameter matters in a given context. We will make
use of the fact that
\[
\Xp(a)^{-1} = \Xp\bigl(-a/(1+pa)\bigr)\vadjust{\goodbreak}
\]
is also an exchangeable design. The function $c_p(a) = -a/(1+pa)$
maps the interval $(-1/p,\infty)$ onto itself, and it holds
$c_p(0)=0$, $c_p(a)\downarrow-1/p$ as $a \uparrow+\infty$, and vice
versa. Exchangeable designs include orthogonal designs for $a=0$, and
they extend to two types of strict collinearities: for $a \uparrow
\infty$ the predictor vectors collapse to a single dimension
$\spanof(\1)$, and for $a \downarrow-1/p$ they collapse to a subspace
$\spanof(\1)^\perp$ of dimension $(p-1)$, where $\1=(1,1,\ldots,1)^T \in
\reals^p$.

As collinearity drives the fracturing of the regression coefficients
into model-dependent quantities $\beta_{j \cdot\rM}$, it is of
interest to analyze $K(\X(a))$ as $\X(a)$ moves from orthogonality at
$a=0$ toward either of the two types of collinearity. Here is what we
find: unguided intuition might suggest that the collapse to rank 1
calls for larger $K(\X)$ than the collapse to rank $p - 1$. This
turns out to be entirely wrong: collapse to rank 1 or rank $p - 1$
has identical effects on $K(\X)$. The reason is duality
(Section~\ref{secdual}): for exchangeable designs, $\X(a)$ collapses
to rank 1 if and only if $\X(a)^* = \X(a)^{-1} = \X(-a/(1+pa))$
collapses to
rank $p - 1$, and vice versa, while $K(\X(a)^{-1}) = K(\X(a))$
according to Corollary~\ref{symdualK}.

We next address the asymptotic behavior of $K=K(\Xp,\alpha)$ for
increasing~$p$. As noted in Section~\ref{secscheffe}, there is a
wide gap between orthogonal designs with $\Korth\sim\sqrt{2 \log
p}$ and the full Scheff\'{e} protection with $\KSch\sim
\sqrt{p}$. The following theorem shows how exchangeable designs fall
into this gap:

%th6.1 #&#
\begin{theorem}\label{thmex}
PoSI constants of exchangeable design matrices $\Xp(a)$ [defined in
(\ref{eqex}) above] have the following limiting behavior:
\[
\lim_{p \rightarrow\infty} \sup_{a \in(-1/p,\infty)} \frac
{K(\Xp(a),\alpha)}{\sqrt{2 \log{p}}} =
2.
\]
\end{theorem}

The proof can be found in Appendix~\ref{appthmex}. The theorem shows
that for exchangeable designs the PoSI constant remains much closer to
the orthogonal case than the Scheff\'{e} case. Thus, for this family
of designs it is possible to improve on the Scheff\'{e} constant by a
considerable margin.

The following detail of geometry for exchangeable designs has a
bearing on their PoSI constants: the angle between pairs of predictor
vectors as a function of $a$ is $\cos(\X_j^{(p)}(a),\X_{j'}^{(p)}(a))
= a(2+pa)/(pa^2+4a+2)$. As the vectors fall into the rank-$(p - 1)$
collinearity at $a \downarrow-1/p$, the cosine becomes
$-1/(2p - 3)$, which converges to zero as $p \rightarrow\infty$.
Thus, as $p \uparrow\infty$, exchangeable designs approach orthogonal
designs even at their most collinear extreme. For further
illustrative materials related to exchangeable designs, see
Figures~2 and 3 in the
online Appendix of the supplementary material [\citet{Beretal},
Section B.7].

%s6.2 #&#
\subsection{\texorpdfstring{Example 2: Where $K(\X)$ is close to the Scheff\'{e} bound.}{Example 2: Where $K(\X)$ is close to the Scheff\'{e}
bound}}
\label{secrootp}

The following is a situation in which the asymptotic upper bound for
$K(\Xp,\alpha)$ is $O(\sqrt{p})$, hence equal to the rate of the\vadjust{\goodbreak}
Scheff\'{e} constant $\KSch(\alpha,p)$. Perhaps surprisingly, it is
sufficient to consider PoSI1 (Section~\ref{secposi1}) whose constant
is dominated by that of full PoSI. Let the PoSI1 predictor of
interest be $\X_p^{(p)}$, so the search is over all models $\rM\ni
p$, but inference is sought only for $\beta_{p \cdot\rM}$. Consider
the following upper triangular $p \times p$ design matrix in canonical
coordinates:
%
%e6.2 #&#
\begin{equation}
\label{eqppcdesign} \Xp(c) = \bigl({\mathbf{e}}_1, {
\mathbf{e}}_2,\ldots, {\mathbf{e}}_{p-1}, \X_p(c)
\bigr),
\end{equation}
where $\X_p(c) = (c, c,\ldots, c, \sqrt{1-(p-1)c^2})^T \in\reals^T$ is
the primary predictor, and the canonical basis vectors ${\mathbf
{e}}_1,\ldots,
{\mathbf{e}}_{p-1} \in\reals^p$ are the controls. The vector $\X
_p(c)$ has
unit length, and hence the parameter $c$ is the correlation between the
primary predictor and the controls. It is constrained to $c^2 <
1/(p-1)$, so $\Xp(c)$ has full rank. For $c^2 = 1/(p-1)$ the primary
predictor $\X_p(c)$ becomes fully collinear with the controls, and it
is on the approach to this boundary where the rate of the following
theorem is attained:

%th6.2 #&#
\begin{theorem} \label{thmrootp} For $\sigma^2$ known, the designs
(\ref{eqppcdesign}) have PoSI\textit{1} constants $K_{p
\cdot}(\Xp(c),\alpha)$ with the following asymptotic rate:
\[
\lim_{p \rightarrow\infty} \sup_{c^2 < 1/(p-1)} \frac{K_{p \cdot
}(\Xp(c),\alpha)}{\sqrt{p}} =
0.6363\ldots.
\]
\end{theorem}

The proof is in Appendix~\ref{appthmrootp}. As $K(\X,\alpha) \ge
K_{j \cdot}(\X,\alpha)$ the theorem provides a lower bound on the rate
of the full PoSI constant. The value 0.6363\ldots\ is not maximal, and we
have indications that the supremum over all designs may exceed~0.78.
Together with the upper bound of Corollary~\ref{corbound}, this would
provide a narrow asymptotic range for worst-case PoSI. Most
importantly, the example shows that for some designs PoSI constants
can be much larger than the $O(1)$ $|t|$-quantiles used in common
practice.

%s6.3 #&#
\subsection{\texorpdfstring{Bounding away from Scheff\'{e}.}{Bounding away from
Scheff\'{e}}}
\label{secupper}

The following is a rough asymptotic upper bound on all PoSI constants
$K(\X,\cM,\alpha)$. It has the Scheff\'{e} rate but with a multiplier
that is strictly less than Scheff\'{e}'s. The bound is loose because
it ignores the rich structure of the sets $\cL(\X,\cM)$
(Section~\ref{secorthgeo}) and only uses their cardinality $|\cL|$
(\mbox{$=$}$p 2^{p-1}$ in the classical case $d = p$ and
\mbox{$\cM= \Mall$}).

%th6.3 #&#
\begin{theorem} \label{thmbound}
Denote by $\cL_d$ arbitrary finite sets of $d$-dimensional
unit vectors, $\cL_d \subset S^{d-1}$, such that $|\cL_d| \le a_d$
where $a_d^{1/d} \rightarrow a (>1)$. Denote by $K(\cL_d,\alpha)$
the $(1-\alpha)$-quantile of $\sup_{\ls\in\cL_d} |\ls^T \Z|$.
Then the following describes an asymptotic worst-case bound for
$K(\cL_d,\alpha)$ and its attainment:
\[
\lim_{d \rightarrow\infty} \sup_{|\cL_d|\le a_d} \frac
{K(\cL_d, \alpha)}{\sqrt{d}} =
\biggl( 1-\frac{1}{a^2} \biggr)^{1/2}.
\]
\end{theorem}

The proof of Theorem~\ref{thmbound} (see
Appendix~\ref{appthmbound}) is an adaptation of Wyner's
(\citeyear{Wyn67}) techniques for sphere packing and sphere\vadjust{\goodbreak}
covering. The worst-case bound ($\le$) is based on a surprisingly
crude Bonferroni-style inequality for caps on spheres. Attainment of
the bound ($\ge$) makes use of the artifice of picking the vectors
$\ls\in\cL$ randomly and independently. Applying the theorem to
PoSI sets $\cL= \cL(\X_{n \times p},\Mall)$ in the classical case
$d = p$, we have $|\cL| = p 2^{p-1} = a_p$, hence
$a_p^{1/p} \rightarrow2$, so the theorem applies with $a = 2$:

%co6.1 #&#
\begin{cor} \label{corbound}
In the classical case $d=p$ a universal asymptotic upper bound
for the PoSI constant $K(\X_{n \times p},\Mall,\alpha)$ is
\[
\lim_{p \rightarrow\infty} \sup_{\X_{n \times p}} \frac
{K(\X_{n \times p},\Mall,\alpha)}{\sqrt{p}}
\le\frac{\sqrt{3}}{2} = 0.866\ldots.
\]
\end{cor}

The corollary shows that the asymptotic rate of the PoSI constant, if
it reaches the Scheff\'{e} rate, will always have a multiplier that is
strictly below that of the Scheff\'{e} constant. We do not know
whether there exist designs for which the bound of the corollary is
attained, but the theorem says the bound is sharp for unstructured
sets $\cL$.

%================================================================

%s7 #&#
\section{\texorpdfstring{Summary and discussion.}{Summary and
discussion}}
\label{secconc}

We investigated the Post-Selection Inference or ``PoSI'' problem for
linear models whereby valid statistical tests and confidence intervals
are sought after variable selection, that is, after selecting a subset
of the predictors in a data-driven way. We adopted a framework that
does \textit{not} assume any of the linear models under consideration to
be correct. We allowed the response vector to be centered at an
arbitrary mean vector but with homoscedastic and Gaussian errors. We
further allowed the full predictor matrix $\X_{n \times p}$ to be
rank-deficient, $d = \rank(\X) < p$, and we also allowed the set
$\cM$ of models $\rM$ under consideration to be largely arbitrary. In
this framework we showed that valid post-selection inference is
possible via simultaneous inference. An important enabling principle
is that submodels have their own regression coefficients; put
differently, $\beta_{j \cdot\rM}$ and $\beta_{j \cdot\rM'}$ are
generally different parameters if $\rM\neq\rM'$. We showed that
simultaneity protection for all parameters $\beta_{j \cdot\rM}$
provides valid post-selection inference. In practice this means
enlarging the constant $t_{1-\alpha/2,r}$ used in conventional
inference to a constant $K(\X_{n \times p},\cM,\alpha,r)$ that
provides simultaneity protection for up to $p 2^{p-1}$
parameters $\beta_{j \cdot\rM}$. We showed that the constant depends
strongly on the predictor matrix $\X$ as the asymptotic bound for
$K(\X,\cM,\alpha,r)$ with $d=\rank(\X)$ ranges between the minimum of
$\sqrt{2 \log d}$ achieved for orthogonal designs on the one hand, and
a large fraction of the Scheff\'{e} bound $\sqrt{d}$ on the other
hand. This wide asymptotic range suggests that computation is
critical for problems with large numbers of predictors. In the
classical case $d = p$ our current computational methods are
feasible up to about $p \approx20$.

We carried out post-selection inference in a limited framework.
Several problems remain open, and many natural extensions are
desirable:
\begin{itemize}
\item Among open problems is the quest for the largest fraction of the
asymptotic Scheff\'{e} rate $\sqrt{d}$ attained by PoSI constants.
So far we know this fraction to be at least 0.6363, but no more than
0.8660\ldots\ in the classical case $d = p$. When the size of models
$|\rM|$ is limited as a function of $p$ (``sparse models''), better
rates can be achieved, and we will report these results elsewhere.
\item Computations for $p>20$ are a challenge. Straight enumeration
of the set of up to $p 2^{p-1}$ linear combinations should be
replaced with heuristic shortcuts that yield practically useful
upper bounds on $K(\X_{n \times p}$, $\cM$, $\alpha,r)$ that are
specific to $\X$ and the set of submodels $\cM$, unlike the 0.8660
fraction of the Scheff\'{e} bound which is universal.
\item Situations to which the PoSI framework should be extended
include generalized linear models, mixed effects models, models
with random predictors, as well as prediction problems. Results
for the last two situations will be reported elsewhere.
\item It would be desirable to devise post-selection inference for
specific selection procedures for cases in which a strict model
selection protocol is being adhered to.
\end{itemize}
$R$ code for computing the PoSI constant for up to $p=20$ can be
obtained from the authors' web pages (a manuscript describing the
computations is available from the authors).

\begin{appendix}\label{app}

%s8 #&#
\section*{\texorpdfstring{Appendix: Proofs}{Appendix: Proofs}}

%s8.1 #&#
\subsection{\texorpdfstring{Proof of Proposition \protect\ref{proporth}.}{Proof of Proposition 5.3}}
\label{appproporth}

(1) The matrix $\X_\rM^* = \X_\rM(\X_\rM^T \X_\rM
)^{-1}$ has the vectors
$\l_{j \cdot\rM}$ as its columns. Thus $\l_{j \cdot\rM} \in
\spanof(\X_j\dvtx  j \in\rM)$. Orthogonality $\l_{j \cdot\rM} \perp
\X_{j'}$ for $j' \neq j$ follows from $\X_\rM^T \X_\rM^* = \I_p$. The
same properties hold for the normalized vectors $\ls_{j \cdot\rM}$.

(2) The vectors $\{\ls_{1 \cdot\{1\}}, \ls_{2 \cdot
\{1,2\}}, \ls_{3 \cdot\{1,2,3\}},\ldots, \ls_{p \cdot
\{1,2,\ldots,p\}} \}$ form a Gram--Schmidt series with normalization,
hence they are an o.n. basis of $\reals^p$.

(3) For $\rM\subset\rM'$, $j \in\rM$, $j' \in\rM'
\setminus
\rM$, we have $\ls_{j \cdot\rM} \perp\ls_{j' \cdot\rM}$ because
they can be embedded in an o.n. basis by first enumerating $\rM$ and
subsequently $\rM' \setminus\rM$, with $j$ being last in the
enumeration of $\rM$ and $j'$ last in the enumeration of $\rM'
\setminus\rM$.

(4) For any $(j_0, \rM_0)$, $j_0 \in\rM_0$, there are
$(p-1) 2^{p-2}$ ways to choose a partner $(j_1, \rM_1)$ such that
either $j_1 \in\rM_1 \subset\rM_0 \setminus j_0$ or $\rM_0 \subset
\rM_1 \setminus j_1$, both\vspace*{1pt} of which result in $\ls_{j_0 \cdot\rM_0}
\perp\ls_{j_1 \cdot\rM_1}$ by the previous part.

%s8.2 #&#
\subsection{\texorpdfstring{Proof of duality: Lemma \protect\ref{lemdual} and
Theorem \protect\ref{thmdual}.}{Proof of duality: Lemma 5.1 and Theorem 5.1}}
\label{applemdual}

The proof relies on a careful analysis of orthogonalities as described
in Proposition~\ref{proporth}, part~(3). In what follows we
write $[\A]$ for the column space of a matrix $\A$, and $[\A]^\perp$
for its orthogonal complement. We show first that, for $\rM\cap
\rM^* = \{j\}$, $\rM\cup\rM^* = \rM_F$, the vectors\vadjust{\goodbreak} $\ls^*_{j
\cdot
\rM^*}$ and $\ls_{j \cdot\rM}$ are in the same one-dimensional
subspace, hence are a multiple of each other. To this end we observe:
%
%e8.1 #&#
%e8.2 #&#
%e8.3 #&#
\begin{eqnarray}
&\displaystyle \ls_{j \cdot\rM} \in [\X_\rM],\qquad \ls_{j \cdot\rM} \in[
\X_{\rM\setminus j}]^\perp,&
\\
&\displaystyle \ls^*_{j \cdot\rM^*} \in\bigl[\X^*_{\rM^*}\bigr],\qquad
\ls^*_{j \cdot\rM^*} \in\bigl[\X^*_{\rM^* \setminus j}\bigr]^\perp,&
\\
\label{dualproof3}
&\displaystyle \bigl[\X^*_{\rM^*}\bigr] = [\X_{\rM\setminus j}]^\perp,\qquad
\bigl[\X^*_{\rM^* \setminus j}\bigr]^\perp= [
\X_{\rM}].&
\end{eqnarray}
The first\vspace*{1pt} two lines state that $\ls_{j \cdot\rM}$ and $\ls^*_{j
\cdot
\rM^*}$ are in the respective column spaces of their models, but
orthogonalized with regard to all other predictors in these models.
The last line, which can also be obtained from the orthogonalities
implied by $\X^T \X^* = \I_p$, establishes that the two vectors fall
in the same one-dimensional subspace,
\[
\ls_{j \cdot\rM} \in[\X_{\rM}] \cap[\X_{\rM\setminus j}]^\perp
= \bigl[\X^*_{\rM^*}\bigr] \cap\bigl[\X^*_{\rM^* \setminus j}
\bigr]^\perp\ni\ls^*_{j
\cdot\rM^*}.
\]
Since they\vspace*{1pt} are normalized, it follows $\ls^*_{j \cdot\rM^*} = \pm
\ls_{j \cdot\rM}$. This result is sufficient to imply all of
Theorem~\ref{thmdual}. The lemma, however, makes a slightly stronger
statement involving lengths which we now prove. In order to express
$\l_{j \cdot\rM}$ and $\l^*_{j \cdot\rM^*}$ according to
(\ref{eqposivecs}), we use $\P_{\rM\setminus j}$ as before and we
write $\P^*_{\rM^* \setminus j}$ for the\vspace*{2pt} analogous projection onto the
space spanned by the columns $\rM^* \setminus j$ of $\X^*$. The
method\vspace*{1pt} of proof is to evaluate $\l_{j \cdot\rM}^T \l^*_{j \cdot
\rM^*}$. The main argument is based on
%
%e8.4 #&#
\begin{equation}
\label{dualproof5} \X_j^T (\I-\P_{\rM\setminus j})
\bigl(\I-\P^*_{\rM^* \setminus j}\bigr) \X_j^* = 1,
\end{equation}
which follows from these facts:
\[
\P_{\rM\setminus j} \P^*_{\rM^* \setminus j} = \0,\qquad \P_{\rM\setminus j}
\X_j^* = \0,\qquad \P^*_{\rM^* \setminus j} \X_j = \0,\qquad
\X_j^T \X_j^* = 1,
\]
which in turn are consequences of (\ref{dualproof3}) and $\X^T \X^* =
\I_p$. We also know from (\ref{eqposivecs}) that
%
%e8.5 #&#
\begin{equation}
\label{dualproof6} \|\l_{j \cdot\rM}\| = 1/\bigl\|(\I-\P_{\rM\setminus j})
\X_j\bigr\|,\qquad \bigl\|\l^*_{j \cdot\rM^*}\bigr\| = 1/\bigl\|\bigl(\I-\P^*_{\rM^* \setminus j}
\bigr)\X^*_j\bigr\|.
\end{equation}
Putting together (\ref{dualproof5}), (\ref{dualproof6}) and
(\ref{eqposivecs}), we obtain
%
%e8.6 #&#
\begin{equation}
\label{dualproof7} \l_{j \cdot\rM}^T \l^*_{j \cdot\rM^*} = \|
\l_{j \cdot\rM
}\|^2 \bigl\|\l^*_{j \cdot\rM^*}\bigr\|^2 > 0.
\end{equation}
Because the two vectors are scalar multiples of each other, we also
know~that
%
%e8.7 #&#
\begin{equation}
\label{dualproof8} \l_{j \cdot\rM}^T \l^*_{j \cdot\rM^*} = \pm\|
\l_{j \cdot
\rM}\| \bigl\|\l^*_{j \cdot\rM^*}\bigr\|.
\end{equation}
Putting together (\ref{dualproof7}) and (\ref{dualproof8}) we conclude
\[
\|\l_{j \cdot\rM}\| \bigl\|\l^*_{j \cdot\rM^*}\bigr\| = 1,\qquad \ls^*_{j \cdot\rM^*} =
\ls_{j \cdot\rM}.
\]
This proves the lemma and the theorem.

%s8.3 #&#
\subsection{\texorpdfstring{Proof of Theorem \protect\ref{thmex}.}{Proof of Theorem 6.1}}
\label{appthmex}

The parameter $a$ in equation (\ref{eqex}) can range from $-1/p$ to
$\infty$, but because of duality there is no loss of generality in
considering only the case in which $a \ge0$, and we do so in the
following. Let $\rM\subset\{1,\ldots,p\}$ and $j \in\rM$.\vadjust{\goodbreak}

Consider first the case $|\rM|=1$, hence $\rM= \{j\}$: we have
$\l_{j\cdot\rM} = \X_j$, the $j$th column of $\X$, and $\ls_{j
\cdot
\rM} = \l_{j \cdot\rM} / \sqrt{pa^2+2a+1}$. For any $\Z\in
\reals^p$ it follows\looseness=1
%
%e8.8 #&#
\begin{equation}
\label{eqbdexs1} \bigl|\ls_{j \cdot\rM}^T \Z\bigr| \le|Z_j| +
\biggl\llvert\frac{1}{\sqrt{p}} \sum_k
Z_k \biggr\rrvert\le\|\Z\|_\infty+ \biggl\llvert
\frac{1}{\sqrt{p}} \sum_k Z_k \biggr
\rrvert.
\end{equation}\looseness=0

Consider next the case $|\rM|>1$, and for notational convenience let
$j=1$ and $\rM=\{1,\ldots,m\}$ where $1 < m \le p$. The following
results can then be applied to arbitrary $\rM$ and $j \in\rM$ by
permuting coordinates. The projection of $\X_1$ on the space spanned
by $\X_2,\ldots, \X_m$ must be of the form
\[
\mathrm{Proj} = \frac{c}{m-1} \sum_{k=2}^m
\X_k = \biggl( ca, \underbrace{ca + \frac{c}{m-1},\ldots,ca +
\frac
{c}{m-1}}_{m-1},\underbrace{ca,\ldots, ca}_{p-m}
\biggr),
\]
where the constant $c$ satisfies $\l_{1 \cdot\rM}=(\X_1 -
\mathrm{Proj}) \bot\mathrm{Proj}$. This follows from symmetry, and no
calculation of projection matrices is needed to verify this. Let
$d=1-c$. Then
%
%e8.9 #&#
\begin{equation}
\label{eqbdexl} (\l_{1 \cdot\rM})_k = \cases{ 1+da &\quad $(k=1)$,
\vspace*{2pt}\cr
\displaystyle - \frac{1-d}{m-1} + da &\quad $(2 \le k \le m)$,
\vspace*{2pt}\cr
da &\quad $(k \ge m+1)$.}
\end{equation}
Some algebra starting from $\l_{1 \cdot\rM}^T \X_2=0$ yields
\[
d = \frac{1/(m-1)}{pa^2+2a + 1/(m-1)}.
\]
The term $da$ is nonnegative, maximal wrt $m$ for $m = 2$, and
thereafter maximal wrt $a$ for $a=1/\sqrt{p}$, whence $\max_{a\ge0,
m\ge2} da = 1/(2(\sqrt{p}+1))$ and finally
%
%e8.10 #&#
\begin{equation}
\label{eqbdexda} 0 \le da < \frac{1}{2\sqrt{p}}.
\end{equation}
This fact will make the term $da$ in (\ref{eqbdexl}) asymptotically
irrelevant. Using \mbox{$\|\l_{1 \cdot\rM}\| \ge1$} and $\bar{\l}_{1
\cdot
\rM}= \l_{1 \cdot\rM}/\|\l_{1 \cdot\rM}\|$ as well as
(\ref{eqbdexl}) and (\ref{eqbdexda}) we obtain
%
%e8.11 #&#
\begin{eqnarray}\label{eqbdexs2}
\bigl|\bar{\l}_{1 \cdot\rM}^T \Z\bigr| &\le& |Z_1| +
\frac{1}{m-1} \sum_{j=2}^m |
Z_j | + \Biggl\llvert\frac{1}{2 \sqrt{p}} \sum
_{j=1}^p Z_j \Biggr\rrvert
\nonumber\\[-8pt]\\[-8pt]
&\le& \|\Z\|_\infty+ \|\Z\|_\infty+ \Biggl\llvert
\frac{1}{2 \sqrt{p}} \sum_{j=1}^p
Z_j \Biggr\rrvert.\nonumber
\end{eqnarray}
Combining (\ref{eqbdexs1}) and (\ref{eqbdexs2}) we obtain
for $\Z\sim\cN(\0,\I_p)$ the following:
\begin{eqnarray*}
\sup_{a \ge0; j,\rM\dvtx  j \in\rM} \bigl|\bar{\l}_{j \cdot\rM}^T \Z\bigr| &\le&
2\|\Z\|_\infty+ \Biggl\llvert\frac{1}{\sqrt{p}} \sum
_{j=1}^p Z_j \Biggr\rrvert
\\
&\le& 2 \sqrt{2 \log p}\bigl(1+o_p(1)\bigr) + O_p(1).
\end{eqnarray*}
This verifies that
%
%e8.12 #&#
\begin{equation}
\label{eqbdexineq} \limsup_{p \rightarrow\infty} \frac{\sup_{a \in
(-1/p,\infty
)}K(\X(a)) }{ \sqrt{2 \log p}} \le2 \qquad\mbox{in probability}.
\end{equation}

It remains to prove that equality holds in (\ref{eqbdexineq}). To
this end let $Z_{(1)}<Z_{(2)}<\cdots<Z_{(p)}$ denote the order
statistics of $Z_1$, $Z_2,\ldots,Z_{p}$. Fix $m$. We have in
probability
\[
\lim_{p \rightarrow\infty} \frac{Z_{(1)}}{\sqrt{2 \log p}} = -1
\quad\mbox{and}\quad \lim
_{p \rightarrow\infty} \frac{Z_{(j)}}{\sqrt{2 \log p}} = 1\qquad \forall j\dvtx  p
- m + 2 \le j \le p.
\]
Note that
\[
\lim_{a \rightarrow\infty} da =0 \quad\mbox{and}\quad \lim_{a \rightarrow\infty}
\|\l_{1 \cdot\rM}\|^2 = 1+(m-1)^{-1}.
\]
For a given $\Z$ we choose $\l_{j^* \cdot\rM^*}$ such that
$j^*=j^*(\Z)$ is the index of $Z_{(1)}$ and $\rM^*=\rM^*(\Z)$ includes
$j^*$ as well as the set of indices of $Z_{(k)}$ for $p-m+2 \le k \le
p$. From (\ref{eqbdexl}) we then obtain in probability
\[
\lim_{p \rightarrow\infty, a \rightarrow\infty} \frac{|\bar{\l}_{j^*
\cdot\rM^*}^T \Z| }{ \sqrt{2 \log p}} \ge\frac{2}{\sqrt{1+(m-1)^{-1}}}.
\]
Choosing $m$ arbitrarily large and combining this with
(\ref{eqbdexineq}) yields the desired conclusion.

%s8.4 #&#
\subsection{\texorpdfstring{Proof of Theorem \protect\ref{thmrootp}.}{Proof of Theorem 6.2}}
\label{appthmrootp}

Recall from (\ref{eqppcdesign}) the designs
\[
\X= \bigl({\mathbf{e}}_1, {\mathbf{e}}_3,\ldots, {
\mathbf{e}}_{p-1}, \X_p(c)\bigr),
\]
where $\X_p(c) = (c, c,\ldots, c, \sqrt{1-(p-1)c^2})^T$ is the primary
predictor. The matrix $\X$ will be treated according to PoSI1
(Section~\ref{secposi1}), and hence we will examine the distribution of
$\max_{\rM: p \in\rM} | \ls_{p \cdot\rM} ^T \Z|$ (assuming
$\sigma^2=1$ known). We determine $\ls_{p \cdot\rM}$ for a fixed
model $\rM$ (\mbox{$\ni$}$p$) with $|\rM|=m$,
\[
\ls_{p \cdot\rM,j} = \cases{\displaystyle  \sqrt{{1-(p-1)c^2 \over1-(m-1)c^2}},
&\quad $j=p$,
\vspace*{2pt}\cr
0, &\quad
$j \in\rM\setminus\{p\}$,
\vspace*{2pt}\cr
\displaystyle {c \over\sqrt{1-(m-1)c^2}}, &\quad $j \in
\rM^c$.}
\]
Therefore,
%
%e8.13 #&#
\begin{equation}
\label{eqrootp1}\quad z_{p \cdot\rM} = \ls_{p \cdot\rM}^T\Z=
\sqrt{{1-(p-1)c^2 \over1-(m-1)c^2}} Z_1 + {c \over\sqrt{1-(m-1)c^2}}
\sum
_{j \in\rM^c} Z_j.
\end{equation}
For fixed $m$ we can explicitly maximize the sum on the right-hand side,
\[
\max_{\rM: |\rM|=m} \biggl\llvert\sum_{j \in\rM^c}
Z_j \biggr\rrvert= \max\Biggl( \sum_{j=1}^{p-m}
Z_{(p-j)}, -\sum_{j=1}^{p-m}
Z_{(j)} \Biggr),
\]
where $Z_{(j)}$ is the $j$th order statistic of $Z_1$, $Z_2,\ldots,
Z_{p-1}$, omitting $Z_p$. We can also explicitly maximize the factor
$c/\sqrt{1-(m-1)c^2}$ in (\ref{eqrootp1}),
\[
\sup_{c^2 < 1/(p-1)} {c \over\sqrt{1-(m-1)c^2}} = \frac{1}{\sqrt{p-m}},
\]
and equality is attained as $c^2 \uparrow1/(p-1)$. Therefore, for
fixed $m$, we can continue from (\ref{eqrootp1}) as follows:
\begin{eqnarray*}
&&
\sup_{c^2 < 1/(p-1)} \max_{|\rM|=m} \frac{\llvert z_{p \cdot
\rM} \rrvert}{\sqrt{p}} \\
&&\qquad=
O_p \biggl( \sqrt{\frac{1}{p}} \biggr)
+ \sqrt{\frac{p}{p-m}} \max\Biggl( \sum_{j=1}^{p-m}
Z_{(p-j)} \frac{1}{p}, {-}\sum_{j=1}^{p-m}
Z_{(j)} \frac{1}{p} \Biggr).
\end{eqnarray*}
The reason for writing the two sums in this manner is that we will
interpret them as approximations to Riemann sums. To this end we
borrow from \citet{Bah66} the following approximations
for $j = 1,\ldots,p - 1$:
\[
Z_{(j)} = \Phi^{-1} \biggl(\frac{j}{p} \biggr) +
O_p\bigl(p^{-1/2}\bigr).
\]
Reparametrizing $m = rp$, the anticipated Riemann approximation is
\[
\int_{r}^1 \Phi^{-1}(x) \,dx = \sum
_{j=1}^{p-m} \Phi^{-1} \biggl(
\frac{p-j}{p} \biggr) \frac{1}{p} + O\bigl(p^{-2}\bigr).
\]
Therefore,
\[
\sum_{j=1}^{p-m} Z_{(p-j)}
\frac{1}{p} = \int_{r}^1
\Phi^{-1}(x) \,dx + O_p\bigl(p^{-1/2}\bigr)
\]
and similarly
\[
-\sum_{j=1}^{p-m} Z_{(j)}
\frac{1}{p} = \int_{r}^1
\Phi^{-1}(x) \,dx + O_p\bigl(p^{-1/2}\bigr).
\]
Summarizing,
\begin{eqnarray*}
&&
\sup_{c} \max_{|\rM|=m} \biggl\llvert
\frac{z_{p \cdot\rM
}}{\sqrt{p}} \biggr\rrvert
\\
&&\qquad= \frac{1}{\sqrt{p-m}} \max\Biggl( \sum_{j=1}^{p-m}
Z_{(p-j)} \frac{1}{p}, - \sum_{j=1}^{p-m}
Z_{(j)} \frac{1}{p} \Biggr) + O_p(\sqrt{1/p})
\\
&&\qquad= \frac{1}{\sqrt{1-r}} \int_{r}^1
\Phi^{-1}(x) \,dx + O_p\bigl(p^{-1/2}\bigr) +
O_p(\sqrt{1/p})
\\
&&\qquad= \frac{1}{\sqrt{1-r}} \phi\bigl(\Phi^{-1}(r)\bigr) + O_p(
\sqrt{1/p}).
\end{eqnarray*}
The function $f(r)=\frac{1}{\sqrt{1-r}} \phi(\Phi^{-1}(r))$ is
maximized at $r^* \approx0.72972$ with\break $f(r^*) \approx0.6363277$.
Therefore,
%
%e8.14 #&#
\begin{equation}
\label{eqK1} \limsup_{p \rightarrow\infty} \sup_{c}
\frac{1}{\sqrt{p}} \max_{\rM} \llvert{z_{p \cdot\rM}}\rrvert
= 0.636\ldots.
\end{equation}
The bound is sharp because it is attained by the models that include
the first or last $m^*=r^*p$ order statistics of $\Z$ when $p
\rightarrow\infty$ and $c^2 \uparrow\frac{1}{p-1}$. From
(\ref{eqK1}) we conclude that $K_{1 \cdot}(\X) \sim0.6363 \sqrt{p}$.

%s8.5 #&#
\subsection{\texorpdfstring{Proof of Theorem \protect\ref{thmbound}.}{Proof of Theorem 6.3}} \label{appthmbound}

We show that if $a_p^{1/p} \rightarrow a (> 1)$, then:
\begin{itemize}
\item we have a uniform asymptotic worst-case bound,
\[
\lim_{p \rightarrow\infty} \sup_{|\cL_p| \le a_p} \max
_{\ls
\in\cL_p} \bigl|\ls^T\Z\bigr|/\sqrt{p} \stackrel{\P} {\le}
\sqrt{1-1/a^2},
\]
\item which is attained when $|\cL_p|=a_p$ and $\ls\in\cL_p$ are
i.i.d. $\Unif(S^{p-1})$ independent of $\Z$,
\[
\lim_{p \rightarrow\infty} \max_{\ls\in\cL_p} \bigl|\ls^T\Z
\bigr|/\sqrt{p} \stackrel{\P} {\ge} \sqrt{1-1/a^2}.
\]
\end{itemize}
These facts imply the assertions about $(1-\alpha)$-quantiles
$K(\cL_p)$ of\break ${\max_{\ls\in\cL_p}} |\ls^T\Z|$ in
Theorem~\ref{thmbound}. We decompose $\Z= R \U$ where $R^2 =
\|\Z\|^2 \sim\chi^2_p$ and $\U= \Z/\|\Z\| \sim\Unif(S^{p-1})$ are
independent. Due to $R/\sqrt{p} \stackrel{\P}{\rightarrow} 1$ it
is sufficient to show the following:
\begin{itemize}
\item uniform asymptotic worst-case bound,
%
%e8.15 #&#
\begin{equation}
\label{eqUbound} \lim_{p \rightarrow\infty} \sup_{|\cL_p| \le a_p}
\max_{\ls
\in\cL_p} \bigl| \ls^T \U\bigr| \stackrel{\P} {\le}
\sqrt{1-1/a^2};
\end{equation}
\item attainment of the bound when $|\cL_p| = a_p$ and $\ls\in\cL
_p$ are i.i.d. $\Unif(S^{p-1})$ independent of $\U$,
%
%e8.16 #&#
\begin{equation}
\label{eqLbound} \lim_{p \rightarrow\infty} \max_{\ls\in\cL_p} \bigl|
\ls^T \U\bigr| \stackrel{\P} {\ge} \sqrt{1-1/a^2}.
\end{equation}
\end{itemize}

To show (\ref{eqUbound}), we upper-bound the noncoverage probability
and show that it converges to zero for $K' > \sqrt{1-1/a^2}$. To this
end we start with a Bonferroni-style bound, as in
\citet{Wyn67},
%
%e8.17 #&#
\begin{eqnarray}\label{eqBonf}
\P\Bigl[ \max_{\ls\in\cL} \bigl| \ls^T \U\bigr| > K'
\Bigr] & = &\P\bigcup_{\ls\in\cL} \bigl[ \bigl|\ls^T
\U\bigr| > K' \bigr]
\nonumber
\\
& \le &\sum_{\ls\in\cL} \P\bigl[ \bigl|\ls^T \U\bigr| >
K'\bigr]
\\
& = &|\cL_p| \P\bigl[ |U| > K'\bigr],\nonumber
\end{eqnarray}
where $U$ is any coordinate of $\U$ or projection of $\U$ onto a unit
vector. We will show that bound (\ref{eqBonf}) converges to
zero. We use the fact that $U^2 \sim \operatorname{Beta}(1/2,(p-1)/2)$, hence
%
%e8.18 #&#
\begin{equation}
\label{eqbeta}\quad \P\bigl[ |U| > K'\bigr] = \frac{1}{\Beta(1/2,(p-1)/2)}
\int
_{K'^2}^1 x^{-1/2} (1-x)^{(p-3)/2} \,dx.
\end{equation}
We bound the Beta function and the integral separately,
\[
\frac{1}{\Beta(1/2,(p-1)/2)} = \frac{\Gamma(p/2)}{\Gamma(1/2)\Gamma
((p-1)/2)} < \sqrt{\frac{(p-1)/2}{\pi}},
\]
where we used $\Gamma(x+1/2)/\Gamma(x) < \sqrt{x}$ (a good
approximation, really) and $\Gamma(1/2)=\sqrt{\pi}$.
\[
\int_{K'^2}^1 x^{-1/2}
(1-x)^{(p-3)/2} \,dx \le\frac{1}{K'} \frac
{1}{(p-1)/2}
\bigl(1-K'^2\bigr)^{(p-1)/2},
\]
where we used $x^{-1/2} \le1/K'$ on the integration interval.
Continuing with the chain of bounds from (\ref{eqBonf}) we have
\[
|\cL_p| \P\bigl[ |U| > K'\bigr] \le\frac{1}{K'}
\biggl(\frac{2}{(p-1)\pi} \biggr)^{1/2} \bigl( |\cL_p|^{1/(p-1)}
\sqrt{1-K'^2} \bigr)^{p-1}.
\]
Since $|\cL_p|^{1/(p-1)} \rightarrow a (>0)$, the right-hand side
converges to zero at geometric speed if $a \sqrt{1-K'^2} < 1$, that
is, if $K' > \sqrt{1-1/a^2}$. This proves (\ref{eqUbound}).

To\vspace*{1pt} show (\ref{eqLbound}), we upper-bound the coverage probability and
show that it converges to zero for $K' < \sqrt{1-1/a^2}$. We make use
of independence of $\ls\in\cL_p$, as in
\citet{Wyn67},
%
%e8.19 #&#
\begin{eqnarray}\label{eqexp}
\P\Bigl[ \max_{\ls\in\cL_p} \bigl| \ls^T \U\bigr| \le
K' \Bigr] & = & \prod_{\ls\in\cL_p} \P\bigl[ \bigl|
\ls^T \U\bigr| \le K' \bigr] = \P\bigl[ | U | \le
K' \bigr]^{|\cL_p|}
\nonumber
\\
& = &\bigl(1 - \P\bigl[ | U | > K' \bigr] \bigr)^{|\cL_p|}
\\
& \le &\exp\bigl( - |\cL_p| \P\bigl[ | U | > K' \bigr]
\bigr).\nonumber
\end{eqnarray}
We will lower-bound the probability $\P[ | U | > K' ]$ recalling
(\ref{eqbeta}) and again deal with the Beta function and the
integral separately,
\[
\frac{1}{\Beta(1/2,(p-1)/2)} = \frac{\Gamma(p/2)}{\Gamma(1/2)\Gamma
((p-1)/2)} > \sqrt{\frac{p/2-3/4}{\pi}},
\]
where we used $\Gamma(x+1)/\Gamma(x+1/2) > \sqrt{x+1/4}$ (again, a
good approximation).
\[
\int_{K'^2}^1 x^{-1/2}
(1-x)^{(p-3)/2} \,dx \ge\frac{1}{(p-1)/2} \bigl(1-K'^2
\bigr)^{(p-1)/2},
\]
where we used $x^{-1/2} \ge1$. Putting it all together we bound the
exponent in~(\ref{eqexp}),
\[
|\cL_p| \P\bigl[ |U| > K'\bigr] \ge\frac{\sqrt{p/2-3/4}}{\sqrt{\pi} (p-1)/2}
\bigl( |\cL_p|^{1/(p-1)} \sqrt{1-K'^2}
\bigr)^{p-1}.
\]
Since $|\cL_p|^{1/(p-1)} \rightarrow a (>0)$, the right-hand side converges
to $+\infty$ at nearly geometric speed if $a \sqrt{1-K'^2} > 1$,
that is, $K' < \sqrt{1-1/a^2}$. This proves~(\ref{eqLbound}).
\end{appendix}

\section*{\texorpdfstring{Acknowledgments.}{Acknowledgments}}

$\!\!$We thank E. Candes, L. Dicker, M. Freiman, E. George, A. Krieger,
M. Low, Z. Ma, E. Pitkin, L. Shepp, N. Sloane, P. Shaman and
M. Traskin for very helpful discussions. The acronym ``SPAR'' is due
to M. Freiman. We are indebted to an anonymous reviewer for extensive
and constructive criticism that influenced the positioning of
this article.

\begin{supplement}%[id=suppA]
\stitle{Supplement to ``Valid post-selection inference''\\}
\slink[doi]{10.1214/12-AOS1077SUPP} %[doi,text={...}] - jei reikia
%suskaldyti doi
\sdatatype{.pdf}
\sfilename{aos1077\_supp.pdf}
\sdescription{The online supplement contains the following
  sections:
  \begin{itemize}[B.7]
  \item[B.1] The Full Model Interpretation of Parameters (as a
    contrast to the sub-model interpretation adopted in this article).
  \item[B.2] ``Omitted Variables Bias'' (which is not bias in the
    sense of this article).
  \item[B.3] Proof of Corollary~\ref{corstrong} (strong error control).
  \item[B.4] Alternative PoSI Guarantees (conditional on selection).
  \item[B.5] PoSI P-Value Adjustment for Model Selection.\vadjust{\goodbreak}
  \item[B.6] The PoSI Process [the PoSI problem in terms of a
    $(j,\mathrm{M})$-indexed process].
  \item[B.7] Figures (illustrating PoSI polytopes and results of a
    simulation for exchangeable designs).
  \end{itemize}}
\end{supplement}

% imsref loaded by lrinkeviciute, 2013-03-18 15:32:36
% imsref loaded by lrinkeviciute, 2013-03-18 15:33:24
% imsref loaded by lrinkeviciute, 2013-03-18 15:36:25
%

\printaddresses


\begin{thebibliography}{38}
% BibTex style file: ims.bst, 2013-01-28
% Default style options (sort=0,type=number).
% Used options (sort=1,type=nameyear).

%b1 #&#
\bibitem[\protect\citeauthoryear{Angrist and Pischke}{2009}]{AngPis09}
%
\begin{bbook}[auto:STB|2013/03/04|13:35:07]
\bauthor{\bsnm{Angrist},~\bfnm{J.~D.}\binits{J.~D.}} \AND
\bauthor{\bsnm{Pischke},~\bfnm{J.~S.}\binits{J.~S.}}
(\byear{2009}).
\btitle{Mostly Harmless Econometrics}.
\bpublisher{Princeton Univ. Press}, \blocation{Princeton}.
\bptok{imsref}%
\end{bbook}
%
\endbibitem

%b2 #&#
\bibitem[\protect\citeauthoryear{Bahadur}{1966}]{Bah66}
%
\begin{barticle}[mr]
\bauthor{\bsnm{Bahadur},~\bfnm{R.~R.}\binits{R.~R.}}
(\byear{1966}).
\btitle{A note on quantiles in large samples}.
\bjournal{Ann. Math. Statist.}
\bvolume{37}
\bpages{577--580}.
\bid{issn={0003-4851}, mr={0189095}}
\bptok{imsref}%
\end{barticle}
%
\endbibitem

%b3 #&#
\bibitem[\protect\citeauthoryear{Berk et~al.}{2013}]{Beretal}
%
\begin{bmisc}[auto:STB|2013/03/04|13:35:07]
\bauthor{\bsnm{Berk},~\bfnm{R.}\binits{R.}},
\bauthor{\bsnm{Brown},~\bfnm{L.}\binits{L.}},
\bauthor{\bsnm{Buja},~\bfnm{A.}\binits{A.}},
\bauthor{\bsnm{Zhang},~\bfnm{K.}\binits{K.}} \AND
\bauthor{\bsnm{Zhao},~\bfnm{L.}\binits{L.}}
(\byear{2013}).
\bhowpublished{Supplement to ``Valid post-selection inference.''
DOI:\doiurl{10.1214/12-AOS1077SUPP}}.
\bptok{imsref}%
\end{bmisc}
%
\endbibitem

%b4 #&#
\bibitem[\protect\citeauthoryear{Brown}{1967}]{Bro67}
%
\begin{barticle}[mr]
\bauthor{\bsnm{Brown},~\bfnm{L.}\binits{L.}}
(\byear{1967}).
\btitle{The conditional level of {S}tudent's {$t$} test}.
\bjournal{Ann. Math. Statist.}
\bvolume{38}
\bpages{1068--1071}.
\bid{issn={0003-4851}, mr={0214210}}
\bptok{imsref}%
\end{barticle}
%
\endbibitem

%b5 #&#
\bibitem[\protect\citeauthoryear{Buehler and Feddersen}{1963}]{BueFed63}
%
\begin{barticle}[mr]
\bauthor{\bsnm{Buehler},~\bfnm{R.~J.}\binits{R.~J.}} \AND
\bauthor{\bsnm{Feddersen},~\bfnm{A.~P.}\binits{A.~P.}}
(\byear{1963}).
\btitle{Note on a conditional property of {S}tudent's {$t$}}.
\bjournal{Ann. Math. Statist.}
\bvolume{34}
\bpages{1098--1100}.
\bid{issn={0003-4851}, mr={0150864}}
\bptok{imsref}%
\end{barticle}
%
\endbibitem

%b6 #&#
\bibitem[\protect\citeauthoryear{Claeskens and Hjort}{2003}]{ClaHjo03}
%
\begin{barticle}[mr]
\bauthor{\bsnm{Claeskens},~\bfnm{Gerda}\binits{G.}} \AND
\bauthor{\bsnm{Hjort},~\bfnm{Nils~Lid}\binits{N.~L.}}
(\byear{2003}).
\btitle{The focused information criterion (with discussion)}.
\bjournal{J.~Amer. Statist. Assoc.}
\bvolume{98}
\bpages{900--945}.
\bid{doi={10.1198/016214503000000819}, issn={0162-1459}, mr={2041482}}
\bptok{imsref}%
\end{barticle}
%
\endbibitem

%b7 #&#
\bibitem[\protect\citeauthoryear{Dijkstra and Veldkamp}{1988}]{DijVel88}
%
\begin{bincollection}[auto:STB|2013/03/04|13:35:07]
\bauthor{\bsnm{Dijkstra},~\bfnm{T.~K.}\binits{T.~K.}} \AND
\bauthor{\bsnm{Veldkamp},~\bfnm{J.~H.}\binits{J.~H.}}
(\byear{1988}).
\btitle{Data-driven selection of regressors and the bootstrap}.
In \bbooktitle{On Model Uncertainty and Its Statistical Implications}
(\beditor{\bfnm{T.~K.}\binits{T.~K.}~\bsnm{Dijkstra}}, ed.)
\bpages{17--38}.
\bpublisher{Springer}, \blocation{Berlin}.
\bptok{imsref}%
\end{bincollection}
%
\endbibitem

%b8 #&#
%%
%(\byear{2010}).
%Available at
%%

%b9 #&#
\bibitem[\protect\citeauthoryear{Hall and Carroll}{1989}]{HalCar89}
%
\begin{barticle}[mr]
\bauthor{\bsnm{Hall},~\bfnm{Peter}\binits{P.}} \AND
\bauthor{\bsnm{Carroll},~\bfnm{R.~J.}\binits{R.~J.}}
(\byear{1989}).
\btitle{Variance function estimation in regression: The effect of estimating
the mean}.
\bjournal{J. R. Stat. Soc. Ser. B Stat. Methodol.}
\bvolume{51}
\bpages{3--14}.
\bid{issn={0035-9246}, mr={0984989}}
\bptok{imsref}%
\end{barticle}
%
\endbibitem

%b10 #&#
\bibitem[\protect\citeauthoryear{Hastie, Tibshirani and
Friedman}{2009}]{HasTibFri09}
%
\begin{bbook}[mr]
\bauthor{\bsnm{Hastie},~\bfnm{Trevor}\binits{T.}},
\bauthor{\bsnm{Tibshirani},~\bfnm{Robert}\binits{R.}} \AND
\bauthor{\bsnm{Friedman},~\bfnm{Jerome}\binits{J.}}
(\byear{2009}).
\btitle{The Elements of Statistical Learning: Data Mining, Inference,
and Prediction},
\bedition{2nd} ed.
\bpublisher{Springer}, \blocation{New York}.
\bid{doi={10.1007/978-0-387-84858-7}, mr={2722294}}
\bptok{imsref}%
\end{bbook}
%
\endbibitem

%b11 #&#
\bibitem[\protect\citeauthoryear{Hjort and Claeskens}{2003}]{HjoCla03}
%
\begin{barticle}[mr]
\bauthor{\bsnm{Hjort},~\bfnm{Nils~Lid}\binits{N.~L.}} \AND
\bauthor{\bsnm{Claeskens},~\bfnm{Gerda}\binits{G.}}
(\byear{2003}).
\btitle{Frequentist model average estimators}.
\bjournal{J. Amer. Statist. Assoc.}
\bvolume{98}
\bpages{879--899}.
\bid{doi={10.1198/016214503000000828}, issn={0162-1459}, mr={2041481}}
\bptok{imsref}%
\end{barticle}
%
\endbibitem

%b12 #&#
\bibitem[\protect\citeauthoryear{Ioannidis}{2005}]{Ioa05}
%
\begin{bmisc}[mr]
\bauthor{\bsnm{Ioannidis},~\bfnm{John P.~A.}\binits{J.~P.~A.}}
(\byear{2005}).
\bhowpublished{Why most published research findings are false.
\textit{PLoS Med.} \textbf{2} e124.
DOI:\doiurl{10.1371/journal.pmed.0020124}}.
\bptok{imsref}%
\end{bmisc}
%
\endbibitem

%b13 #&#
\bibitem[\protect\citeauthoryear{Kabaila}{1998}]{Kab98}
%
\begin{barticle}[mr]
\bauthor{\bsnm{Kabaila},~\bfnm{Paul}\binits{P.}}
(\byear{1998}).
\btitle{Valid confidence intervals in regression after variable selection}.
\bjournal{Econometric Theory}
\bvolume{14}
\bpages{463--482}.
\bid{doi={10.1017/S0266466698144031}, issn={0266-4666}, mr={1650037}}
\bptok{imsref}%
\end{barticle}
%
\endbibitem

%b14 #&#
\bibitem[\protect\citeauthoryear{Kabaila}{2009}]{Kab09}
%
\begin{barticle}[auto:STB|2013/03/04|13:35:07]
\bauthor{\bsnm{Kabaila},~\bfnm{P.}\binits{P.}}
(\byear{2009}).
\btitle{The coverage properties of confidence regions after model selection}.
\bjournal{International Statistical Review}
\bvolume{77}
\bpages{405--414}.
\bptok{imsref}%
\end{barticle}
%
\endbibitem

%b15 #&#
\bibitem[\protect\citeauthoryear{Kabaila and Leeb}{2006}]{KabLee06}
%
\begin{barticle}[mr]
\bauthor{\bsnm{Kabaila},~\bfnm{Paul}\binits{P.}} \AND
\bauthor{\bsnm{Leeb},~\bfnm{Hannes}\binits{H.}}
(\byear{2006}).
\btitle{On the large-sample minimal coverage probability of confidence
intervals after model selection}.
\bjournal{J. Amer. Statist. Assoc.}
\bvolume{101}
\bpages{619--629}.
\bid{doi={10.1198/016214505000001140}, issn={0162-1459}, mr={2256178}}
\bptok{imsref}%
\end{barticle}
%
\endbibitem

%b16 #&#
\bibitem[\protect\citeauthoryear{Leeb}{2006}]{Lee06}
%
\begin{bincollection}[mr]
\bauthor{\bsnm{Leeb},~\bfnm{Hannes}\binits{H.}}
(\byear{2006}).
\btitle{The distribution of a linear predictor after model selection:
Unconditional finite-sample distributions and asymptotic approximations}.
In \bbooktitle{Optimality}.
\bseries{Institute of Mathematical Statistics Lecture
Notes---Monograph Series}
\bvolume{49}
\bpages{291--311}.
\bpublisher{IMS}, \blocation{Beachwood, OH}.
\bid{doi={10.1214/074921706000000518}, mr={2338549}}
\bptok{imsref}%
\end{bincollection}
%
\endbibitem

%b17 #&#
\bibitem[\protect\citeauthoryear{Leeb and P{\"o}tscher}{2003}]{LeePot03}
%
\begin{barticle}[mr]
\bauthor{\bsnm{Leeb},~\bfnm{Hannes}\binits{H.}} \AND
\bauthor{\bsnm{P{\"o}tscher},~\bfnm{Benedikt~M.}\binits{B.~M.}}
(\byear{2003}).
\btitle{The finite-sample distribution of post-model-selection
estimators and
uniform versus nonuniform approximations}.
\bjournal{Econometric Theory}
\bvolume{19}
\bpages{100--142}.
\bid{doi={10.1017/S0266466603191050}, issn={0266-4666}, mr={1965844}}
\bptok{imsref}%
\end{barticle}
%
\endbibitem

%b18 #&#
\bibitem[\protect\citeauthoryear{Leeb and P{\"o}tscher}{2005}]{LeePot05}
%
\begin{barticle}[mr]
\bauthor{\bsnm{Leeb},~\bfnm{Hannes}\binits{H.}} \AND
\bauthor{\bsnm{P{\"o}tscher},~\bfnm{Benedikt~M.}\binits{B.~M.}}
(\byear{2005}).
\btitle{Model selection and inference: Facts and fiction}.
\bjournal{Econometric Theory}
\bvolume{21}
\bpages{21--59}.
\bid{doi={10.1017/S0266466605050036}, issn={0266-4666}, mr={2153856}}
\bptok{imsref}%
\end{barticle}
%
\endbibitem

%b19 #&#
\bibitem[\protect\citeauthoryear{Leeb and P{\"o}tscher}{2006a}]{LeePot06N1}
%
\begin{barticle}[mr]
\bauthor{\bsnm{Leeb},~\bfnm{Hannes}\binits{H.}} \AND
\bauthor{\bsnm{P{\"o}tscher},~\bfnm{Benedikt~M.}\binits{B.~M.}}
(\byear{2006}a).
\btitle{Performance limits for estimators of the risk or distribution of
shrinkage-type estimators, and some general lower risk-bound results}.
\bjournal{Econometric Theory}
\bvolume{22}
\bpages{69--97}.
\bid{doi={10.1017/S0266466606060038}, issn={0266-4666}, mr={2212693}}
\bptok{imsref}%
\end{barticle}
%
\endbibitem

%b20 #&#
\bibitem[\protect\citeauthoryear{Leeb and P{\"o}tscher}{2006b}]{LeePot06N2}
%
\begin{barticle}[mr]
\bauthor{\bsnm{Leeb},~\bfnm{Hannes}\binits{H.}} \AND
\bauthor{\bsnm{P{\"o}tscher},~\bfnm{Benedikt~M.}\binits{B.~M.}}
(\byear{2006}b).
\btitle{Can one estimate the conditional distribution of post-model-selection
estimators?}
\bjournal{Ann. Statist.}
\bvolume{34}
\bpages{2554--2591}.
\bid{doi={10.1214/009053606000000821}, issn={0090-5364}, mr={2291510}}
\bptok{imsref}%
\end{barticle}
%
\endbibitem

%b21 #&#
\bibitem[\protect\citeauthoryear{Leeb and P{\"o}tscher}{2008a}]{LeePot08N1}
%
\begin{bincollection}[auto:STB|2013/03/04|13:35:07]
\bauthor{\bsnm{Leeb},~\bfnm{H.}\binits{H.}} \AND
\bauthor{\bsnm{P{\"o}tscher},~\bfnm{B.~M.}\binits{B.~M.}}
(\byear{2008}a).
\btitle{Model selection}.
In \bbooktitle{The Handbook of Financial Time Series}
(\beditor{\bfnm{T.~G.}\binits{T.~G.}~\bsnm{Anderson}},
\beditor{\bfnm{R.~A.}\binits{R.~A.}~\bsnm{Davis}},
\beditor{\bfnm{J.~P.}\binits{J.~P.}~\bsnm{Kreiss}} \AND
\beditor{\bfnm{T.}\binits{T.}~\bsnm{Mikosch}}, eds.)
\bpages{785--821}.
\bpublisher{Springer}, \blocation{New York}.
\bptok{imsref}%
\end{bincollection}
%
\endbibitem

%b22 #&#
\bibitem[\protect\citeauthoryear{Leeb and P{\"o}tscher}{2008b}]{LeePot08N2}
%
\begin{barticle}[mr]
\bauthor{\bsnm{Leeb},~\bfnm{Hannes}\binits{H.}} \AND
\bauthor{\bsnm{P{\"o}tscher},~\bfnm{Benedikt~M.}\binits{B.~M.}}
(\byear{2008}b).
\btitle{Can one estimate the unconditional distribution of post-model-selection
estimators?}
\bjournal{Econometric Theory}
\bvolume{24}
\bpages{338--376}.
\bid{doi={10.1017/S0266466608080158}, issn={0266-4666}, mr={2422862}}
\bptok{imsref}%
\end{barticle}
%
\endbibitem

%b23 #&#
\bibitem[\protect\citeauthoryear{Leeb and P{\"o}tscher}{2008c}]{LeePot08N3}
%
\begin{barticle}[mr]
\bauthor{\bsnm{Leeb},~\bfnm{Hannes}\binits{H.}} \AND
\bauthor{\bsnm{P{\"o}tscher},~\bfnm{Benedikt~M.}\binits{B.~M.}}
(\byear{2008}c).
\btitle{Sparse estimators and the oracle property, or the return of {H}odges'
estimator}.
\bjournal{J. Econometrics}
\bvolume{142}
\bpages{201--211}.
\bid{doi={10.1016/j.jeconom.2007.05.017}, issn={0304-4076}, mr={2394290}}
\bptok{imsref}%
\end{barticle}
%
\endbibitem

%b24 #&#
\bibitem[\protect\citeauthoryear{Moore and McCabe}{2003}]{MooMcC03}
%
\begin{bbook}[auto:STB|2013/03/04|13:35:07]
\bauthor{\bsnm{Moore},~\bfnm{D.~S.}\binits{D.~S.}} \AND
\bauthor{\bsnm{McCabe},~\bfnm{G.~P.}\binits{G.~P.}}
(\byear{2003}).
\btitle{Introduction to the Practice of Statistics},
\bedition{4th} ed.
\bpublisher{Freeman}, \blocation{New York}.
\bptok{imsref}%
\end{bbook}
%
\endbibitem

%b25 #&#
\bibitem[\protect\citeauthoryear{Olshen}{1973}]{Ols73}
%
\begin{barticle}[mr]
\bauthor{\bsnm{Olshen},~\bfnm{Richard~A.}\binits{R.~A.}}
(\byear{1973}).
\btitle{The conditional level of the {$F$}-test}.
\bjournal{J. Amer. Statist. Assoc.}
\bvolume{68}
\bpages{692--698}.
\bid{issn={0162-1459}, mr={0359198}}
\bptok{imsref}%
\end{barticle}
%
\endbibitem

%b26 #&#
\bibitem[\protect\citeauthoryear{P{\"o}tscher}{1991}]{Pot91}
%
\begin{barticle}[mr]
\bauthor{\bsnm{P{\"o}tscher},~\bfnm{B.~M.}\binits{B.~M.}}
(\byear{1991}).
\btitle{Effects of model selection on inference}.
\bjournal{Econometric Theory}
\bvolume{7}
\bpages{163--185}.
\bid{doi={10.1017/S0266466600004382}, issn={0266-4666}, mr={1128410}}
\bptok{imsref}%
\end{barticle}
%
\endbibitem

%b27 #&#
\bibitem[\protect\citeauthoryear{P{\"o}tscher}{2006}]{Pot06}
%
\begin{bincollection}[mr]
\bauthor{\bsnm{P{\"o}tscher},~\bfnm{Benedikt~M.}\binits{B.~M.}}
(\byear{2006}).
\btitle{The distribution of model averaging estimators and an impossibility
result regarding its estimation}.
In \bbooktitle{Time Series and Related Topics}.
\bseries{Institute of Mathematical Statistics Lecture
Notes---Monograph Series}
\bvolume{52}
\bpages{113--129}.
\bpublisher{IMS}, \blocation{Beachwood, OH}.
\bid{doi={10.1214/074921706000000987}, mr={2427842}}
\bptok{imsref}%
\end{bincollection}
%
\endbibitem

%b28 #&#
\bibitem[\protect\citeauthoryear{P{\"o}tscher and Leeb}{2009}]{PotLee09}
%
\begin{barticle}[mr]
\bauthor{\bsnm{P{\"o}tscher},~\bfnm{Benedikt~M.}\binits{B.~M.}}
\AND
\bauthor{\bsnm{Leeb},~\bfnm{Hannes}\binits{H.}}
(\byear{2009}).
\btitle{On the distribution of penalized maximum likelihood
estimators: The
{LASSO}, {SCAD}, and thresholding}.
\bjournal{J. Multivariate Anal.}
\bvolume{100}
\bpages{2065--2082}.
\bid{doi={10.1016/j.jmva.2009.06.010}, issn={0047-259X}, mr={2543087}}
\bptok{imsref}%
\end{barticle}
%
\endbibitem

%b29 #&#
\bibitem[\protect\citeauthoryear{P{\"o}tscher and Schneider}{2009}]{PotSch09}
%
\begin{barticle}[mr]
\bauthor{\bsnm{P{\"o}tscher},~\bfnm{Benedikt~M.}\binits{B.~M.}}
\AND
\bauthor{\bsnm{Schneider},~\bfnm{Ulrike}\binits{U.}}
(\byear{2009}).
\btitle{On the distribution of the adaptive {LASSO} estimator}.
\bjournal{J. Statist. Plann. Inference}
\bvolume{139}
\bpages{2775--2790}.
\bid{doi={10.1016/j.jspi.2009.01.003}, issn={0378-3758}, mr={2523666}}
\bptok{imsref}%
\end{barticle}
%
\endbibitem

%b30 #&#
\bibitem[\protect\citeauthoryear{P{\"o}tscher and Schneider}{2010}]{PotSch10}
%
\begin{barticle}[mr]
\bauthor{\bsnm{P{\"o}tscher},~\bfnm{Benedikt~M.}\binits{B.~M.}}
\AND
\bauthor{\bsnm{Schneider},~\bfnm{Ulrike}\binits{U.}}
(\byear{2010}).
\btitle{Confidence sets based on penalized maximum likelihood
estimators in
{G}aussian regression}.
\bjournal{Electron. J. Stat.}
\bvolume{4}
\bpages{334--360}.
\bid{doi={10.1214/09-EJS523}, issn={1935-7524}, mr={2645488}}
\bptok{imsref}%
\end{barticle}
%
\endbibitem

%b31 #&#
\bibitem[\protect\citeauthoryear{P{\"o}tscher and Schneider}{2011}]{PotSch11}
%
\begin{barticle}[mr]
\bauthor{\bsnm{P{\"o}tscher},~\bfnm{Benedikt~M.}\binits{B.~M.}}
\AND
\bauthor{\bsnm{Schneider},~\bfnm{Ulrike}\binits{U.}}
(\byear{2011}).
\btitle{Distributional results for thresholding estimators in high-dimensional
{G}aussian regression models}.
\bjournal{Electron. J. Stat.}
\bvolume{5}
\bpages{1876--1934}.
\bid{doi={10.1214/11-EJS659}, issn={1935-7524}, mr={2970179}}
\bptok{imsref}%
\end{barticle}
%
\endbibitem

%b32 #&#
%%
%(\byear{2006}).
%Applications}.
%%

%b33 #&#
\bibitem[\protect\citeauthoryear{Scheff{\'e}}{1959}]{Sch59}
%
\begin{bbook}[mr]
\bauthor{\bsnm{Scheff{\'e}},~\bfnm{Henry}\binits{H.}}
(\byear{1959}).
\btitle{The Analysis of Variance}.
\bpublisher{Wiley}, \blocation{New York}.
\bid{mr={0116429}}
\bptok{imsref}%
\end{bbook}
%
\endbibitem

%b34 #&#
\bibitem[\protect\citeauthoryear{Sen}{1979}]{Sen79}
%
\begin{barticle}[mr]
\bauthor{\bsnm{Sen},~\bfnm{Pranab~Kumar}\binits{P.~K.}}
(\byear{1979}).
\btitle{Asymptotic properties of maximum likelihood estimators based on
conditional specification}.
\bjournal{Ann. Statist.}
\bvolume{7}
\bpages{1019--1033}.
\bid{issn={0090-5364}, mr={0536504}}
\bptok{imsref}%
\end{barticle}
%
\endbibitem

%b35 #&#
\bibitem[\protect\citeauthoryear{Sen and Saleh}{1987}]{SenSal87}
%
\begin{barticle}[mr]
\bauthor{\bsnm{Sen},~\bfnm{Pranab~Kumar}\binits{P.~K.}} \AND
\bauthor{\bsnm{Saleh},~\bfnm{A.~K. M.~Ehsanes}\binits{A.~K. M.~E.}}
(\byear{1987}).
\btitle{On preliminary test and shrinkage {$M$}-esti\-mation in linear models}.
\bjournal{Ann. Statist.}
\bvolume{15}
\bpages{1580--1592}.
\bid{doi={10.1214/aos/1176350611}, issn={0090-5364}, mr={0913575}}
\bptok{imsref}%
\end{barticle}
%
\endbibitem

%b36 #&#
%%
%(\byear{2009}).
%%\bnote{Revised and extended from the 2004 French original, Translated
%%by
%% Vladimir Zaiats}.
%%

%b37 #&#
\bibitem[\protect\citeauthoryear{Wyner}{1967}]{Wyn67}
%
\begin{barticle}[mr]
\bauthor{\bsnm{Wyner},~\bfnm{A.~D.}\binits{A.~D.}}
(\byear{1967}).
\btitle{Random packings and coverings of the unit {$n$}-sphere}.
\bjournal{Bell System Tech. J.}
\bvolume{46}
\bpages{2111--2118}.
\bid{issn={0005-8580}, mr={0223979}}
\bptok{imsref}%
\end{barticle}
%
\endbibitem

\end{thebibliography}
\end{document}